\documentclass[12pt,reqno]{amsart}

\usepackage{fullpage}

\usepackage{times}
\usepackage[T1]{fontenc}
\usepackage{mathrsfs}
\usepackage{latexsym}
\usepackage[dvips]{graphics}
\usepackage[dvips]{graphicx}
\usepackage{epsfig}
\usepackage{amsmath,amsfonts,amsthm,amssymb,amscd}
\input amssym.def
\input amssym.tex
\usepackage{color}
\usepackage{nicefrac}

\usepackage{verbatim}

\usepackage{epstopdf}

\usepackage{tikz}
\usetikzlibrary{arrows}
\tikzstyle{block}=[draw opacity=0.7,line width=1.4cm]

\usepackage[shortalphabetic,msc-links]{amsrefs}
\newcommand{\pp}{\ensuremath{\mathbb{P}}}

\newcommand\be{\begin{equation}}
\newcommand\ee{\end{equation}}
\newcommand\bea{\begin{eqnarray}}
\newcommand\eea{\end{eqnarray}}
\newcommand\nbea{\begin{eqnarray*}}
\newcommand\neea{\end{eqnarray*}}
\newcommand\bi{\begin{itemize}}
\newcommand\ei{\end{itemize}}
\newcommand\ben{\begin{enumerate}}
\newcommand\een{\end{enumerate}}

\newtheorem{thm}{Theorem}[section]
\newtheorem{conj}[thm]{Conjecture}
\newtheorem{cor}[thm]{Corollary}
\newtheorem{lem}[thm]{Lemma}
\newtheorem{prop}[thm]{Proposition}

\newtheorem{rek}[thm]{Remark}

\newcommand{\E}{\ensuremath{\mathbb{E}}}

\numberwithin{equation}{section}

\DeclareMathOperator{\Var}{Var}

\newcommand{\bal}{\begin{align}}
\newcommand{\eal}{\end{align}}

\newcommand{\NN}{\ensuremath{\mathbb N}}
\newcommand{\ZZ}{\ensuremath{\mathbb Z}}
\newcommand{\Prob}[1]{{\mathbb P}\left( #1 \right)}
\newcommand{\Expect}[1]{{\mathbb E}\left[ #1 \right]}

\begin{document}

\title[Distribution of Missing Sums in Sumsets]{Distribution of Missing Sums in Sumsets}

\author{Oleg Lazarev}\email{olazarev@Princeton.edu}
\address{Department of Mathematics, Princeton University, Princeton, NJ 08544}

\author{Steven J. Miller}\email{sjm1@williams.edu, Steven.Miller.MC.96@aya.yale.edu}
\address{Department of Mathematics and Statistics, Williams College,
Williamstown, MA 01267}

\author{Kevin O'Bryant}\email{kevin@member.ams.org}
\address{Department of Mathematics, CUNY, The College of Staten Island and the Graduate Center,
Staten Island, NY 10314}



\date{\today}

\thanks{We thank the participants of the SMALL 2011 REU at Williams College for many enlightening conversations, and the referee for many helpful comments on an earlier draft. The first named author was supported by NSF grants DMS0850577 and Williams College; the second named author was partially supported by NSF grant DMS0970067. This research was supported, in part, under National Science Foundation Grants CNS-0958379 and CNS-0855217 and the City University of New York High Performance Computing Center.}

\begin{abstract}
For any finite set of integers $X$, define its sumset $X+X$ to be $\{x+y: x, y \in X\}$. In a recent paper, Martin and O'Bryant investigated the distribution of $|A+A|$ given the uniform distribution on subsets $A\subseteq \{0, 1, \dots, n-1\}$. They also conjectured the existence of a limiting distribution for $|A+A|$ and showed that the expectation of $|A+A|$ is $2n - 11 + O((3/4)^{n/2})$. Zhao proved that the limits $m(k):=\lim_{n\to\infty} \Prob{2n-1-|A+A|=k}$ exist, and that $\sum_{k\geq0} m(k)=1$.

We continue this program and give exponentially decaying upper and lower bounds on $m(k)$, and sharp bounds on $m(k)$ for small $k$. Surprisingly, the distribution is at least bimodal; sumsets have an unexpected bias against missing exactly 7 sums. The proof of the latter is by reduction to questions on the distribution of related random variables, with large scale numerical computations a key ingredient in the analysis. We also derive an explicit formula for the variance of $|A+A|$ in terms of Fibonacci numbers, finding $\Var(|A+A|) \approx 35.9658$. New difficulties arise in the form of weak dependence between events of the form $\{x\in A+A\}$, $\{y\in A+A\}$. We surmount these obstructions by translating the problem to graph theory. This approach also yields good bounds on the probability for $A+A$ missing a consecutive block of length $k$.

%For $X \subset \N$, its sumset $X+X$ is $\{x+y: x, y \in X\}$. Martin \& O'Bryant investigated the distribution of $|A+A|$ given the uniform distribution on subsets $A\subseteq \{0, 1, \dots, n-1\}$. They conjectured the limiting distribution exists and proved $\E[|A+A|]$ $=$ $2n - 11 + O((3/4)^{n/2})$; later Zhao proved each $m(k):=\lim_{n\to\infty} \Prob{2n-1-|A+A|=k}$ exists. We derive exponentially decaying upper and lower bounds on $m(k)$. Surprisingly, the distribution is at least bimodal; sumsets have an unexpected bias against missing exactly 7 sums. We express ${\rm Var}\left(|A+A|\right)$ in terms of Fibonacci numbers. Difficulties arise from dependence between events $\{x\in A+A\}$, $\{y\in A+A\}$. We surmount these by translating the problem to graph theory.
\end{abstract}

%These bounds are based on bounds on probabilities like $\pp(k+ a_1, \dots, \mbox{ and } k + a_m \not \in A+A)$, which we prove are roughly exponential in $k$ for fixed $a_1, \dots, a_m$ by applying a modified version of Fekete's Lemma.

\maketitle

\tableofcontents

%\clearpage

%%%%%%%%%%%%%%%%%%%%%%%%%%%%%%%%%%%%%%%%%%%%%%%%%%%%%%%%%%%%%%%%%%%%%%%%%%%%%%%%%%
%%%%%%%%%%%%%%%%%%%%%%%%%%%%%%%%%%%%%%%%%%%%%%%%%%%%%%%%%%%%%%%%%%%%%%%%%%%%%%%%%%

\section{Introduction}
The central object of additive number theory \cites{Na,TV} is the sumset $X+X$ of a set $X$ of integers:
\begin{equation}
   X+X\ :=\ \{x_1 + x_2 : x_1, x_2 \in X \}.
\end{equation}
Typically, the theory is concerned with extremal behavior, such as the structure of finite $X$ when $|X+X|/|X|$ is nearly minimal (Freiman's Theorem), or the possible densities of $X$ when $|X+X|/\tbinom{|X|+1}{2}$ is maximized (Sidon Sets). See \cites{Na,Ru} for surveys and \cites{Fr,Ji} for examples.

Here we focus on {\em typical} behavior: for a randomly chosen set $X$ of integers, what is the expected value and variance of $|X+X|$? The answer of course depends on how $X$ is chosen, and we focus our attention on sets taken uniformly from the $2^n$ subsets of $[0,n-1]$; we denote intervals of integers as $[a,b]:=\{x\in \ZZ : a\leq x \leq b\}$ and such a random set as $A$.
In \S\ref{subsec: intro divot} and \S\ref{subsec: connecting the distribs} we discuss some variations on the manner of choosing a random set of natural numbers.

Other authors have considered aspects of typical behavior of sumsets. When Erd\H{o}s and R\'{e}nyi \cite{ErdosRenyi} first applied the probabilistic method to number theory, they observed that with probability 1, a uniformly random subset $C$ of $\NN$ will have $C+C=\NN\setminus F$ for some {\em finite} set $F$, but made no effort to explore $F$ further. The present work concerns itself with properties of the set $$F_n:=[0,2n-2]\setminus(A+A),$$ with $A$ as above. We prove the existence of $$\lim_{n\to\infty} \Expect{|F_n|^r}$$ for every $r\geq 1$, give upper and lower bounds on $$\Prob{|F_n|=k}$$ for small $k$, large $n$, and also as $k\to\infty$, and also bound $$\Prob{\{a_1, a_2, \dots, a_k\} \subseteq F_n}.$$ Our work is usually quantitatively effective, and we report numerical estimates throughout.

The key obstacle to finding the limiting distribution of $|F_n|$ is the dependence between different elements occurring or not occurring in $A+A$. For example, $3 \not \in A+A$ and $7 \not \in A+A$ are dependent events since both are affected by whether $2 \in A$. We develop a graph theoretic framework which makes it much easier to analyze the dependence between such events and to develop bounds that incorporate the dependence. It is possible to avoid this framework, but doing so makes both notation and the underlying issues less clear.

Graph theory has been used in additive number theory before. For example, Pl\"unnecke (see the description in \cite{Ru}) uses graph theory to estimate the size of $k$-fold sumsets in terms of $|A|$ and $|A+A|$, Alon and Erd\"{o}s \cite{AE} use hypergraphs to study Sidon sets, and Gilbert \cite{Gi} on the Erd{\H{o}}s-Turan conjecture. Our use of graph theory seems to be different from these
as we investigate the size of $A+A$ for typical $A$, without reference to the size of $A$ itself.

The next subsection of this introduction sets up our notation and states our main results. The last two subsections provide more motivation and indicate the nature of our proofs and computations. In \S\ref{sec: graph}, we develop a graph theoretic framework for handling the dependencies between events like $\{a_1\in F\}$ and $\{a_2\in F\}$. In \S\ref{sec: variance}, we find an explicit formula for the limit of the variance of $|F|$ and prove Theorem \ref{thm: variance}, stated below. In \S\ref{sec:  bounds}, we prove the exponential bounds for Theorem \ref{thm: bounds}. In \S\ref{sec: configuration}, we find the probability of missing certain configurations and prove Theorem \ref{thm: configuration}, while in \S\ref{sec: consecutive} we discuss consecutive missing elements and prove Theorem \ref{thm: consecutive} and Theorem \ref{thm: lambdabounds}. We return to the problem of explicit bounds on $\Prob{|F_n|=k}$ for small $k$ and the existence of a limiting distribution for $|F_n|$ in \S\ref{sec:bimodality}.
Finally in \S\ref{sec: future}, we discuss some problems for future research and how the graph theoretic framework may be applied to such problems.

\begin{rek} Many of the questions in this paper grew out of studying the difference in size between the sumset $A+A$ and the difference set $A-A$. As addition is commutative and subtraction is not, it is natural to expect the difference set of a typical $A$ drawn uniformly from $\{0,1, \dots, n\}$ to be larger than the sumset. Though numerical exploration and heuristics suggested that almost all sets should give rise to more differences, Martin and O'Bryant \cite{MO} proved that a small but positive percentage are sum-dominant. The percentage is quite small, around $4.5 \cdot 10^{-4}$ \cite{Zh}.  Understanding the structure of $A+A$, in particular when and what sums are missing, has motivated much of the theoretical and numerical work in the field. For other directions, see \cite{HM1} for results on non-uniform models or \cite{ILMZ} for multiple comparisons and summands.
\end{rek}

\subsection{Terminology and Theorems}\label{subsec: intro terminology}

The main characteristic of $A+A$ is that it is almost full. Martin and O'Bryant \cite{MO} proved that
\begin{equation}\label{thm: martinobryant}
   \Expect{|A+A|}\ =\ 2n - 1- 10 + O\left((3/4)^{n/2}\right).
\end{equation}
Since typical sumsets are almost full, it is more natural to investigate the number of missing sums, which is why we write the above as $2n-1$ minus 10. As noted in \cite{MO}, sumsets are almost full because middle elements have many representations as a sum of two elements of $[0,n-1]$; each $i \in [0, 2n-2]$ has roughly $n/4 - |n-i|/4$ representations.

\begin{figure}[h]
\begin{center}
\includegraphics[scale=.80]{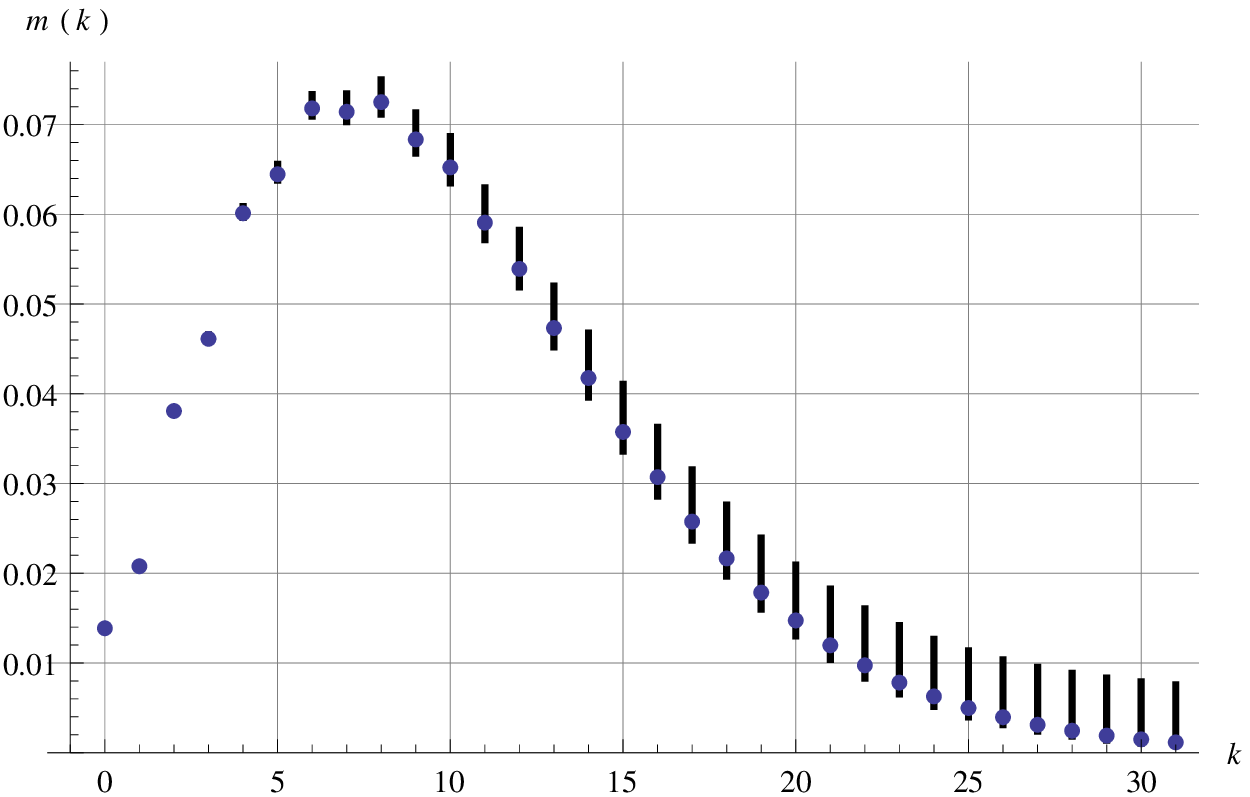}
\end{center}
\caption{Experimental values of $m(k)$, with vertical bars depicting the values allowed by our rigorous bounds. In most cases, the allowed interval is smaller than the dot indicating the experimental value. The data comes from generating $2^{28}$ sets uniformly forced to contain 0 from $[0, 256)$; see \S\ref{subsec: connecting the distribs} for details of the calculation.\label{figure:x_k Experimental graph}}
\end{figure}

We set
   \begin{align}
      M_{[0,n-1]} &\ :=\ \left| [0,2n-2] \setminus (A+A)\right| = 2n-1-|A+A|, \nonumber\\
      m_n(k) &\ :=\ \Prob{M_{[0,n-1]}=k}, \nonumber\\
      m(k) &\ :=\ \lim_{n\to\infty} m_n(k).
   \end{align}
A special case of Zhao's theorem~\cite{Zh} is that $m(k)$ is well-defined, strictly positive, and that $\sum_{k=0}^\infty m(k) =1$, so that we can think of $m(k)$ as defining a distribution on $\NN$. Thus, we can speak of ``the probability that a large finite set $X$ has a sumset that misses exactly 17 elements'' and mean something sensible. Zhao's work is numerically impractical and did not give reasonable upper bounds on $m(k)$; we do that in \S\ref{sec:bimodality}, where we also reprove Zhao's results in this easier setting. See Figure~\ref{figure:x_k Experimental graph} for the experimental estimates and rigorous bounds on $m(k)$ for $0\leq k < 32$.

The result \eqref{thm: martinobryant} above implies that
   \[
      \lim_{n\to\infty} \Expect{M_{[0,n-1]}} = 10.
   \]
Equivalently, in light of Zhao's work, $\sum_{k=0}^\infty k m(k) = 10$. To this, we add the following results. Let $\phi:=(1+\sqrt5)/2$, the golden ratio.

\begin{thm}\label{thm: bounds}
Let $n>5k$. Then
\begin{equation}
2^{-k/2}\ \ll\ m_n(k) \ll\ (\phi/2)^k,
\end{equation}
where the implied constants are independent of $k$ and $n$.
\end{thm}

Note that $2^{-1/2}\approx 0.707$ and $\phi/2\approx 0.809$, so that bounds provided by Theorem~\ref{thm: bounds} are reasonably close. We suspect, based on numerical data, that the following conjecture represents the truth of the matter, and perhaps even $\lambda = \sqrt{\phi-1}$.

\begin{conj}\label{conj: main}
There exists $\lambda$ such that for any $\epsilon > 0$,
\begin{equation}\label{eqn: conj}
(\lambda - \epsilon)^k \ll_\epsilon m(k) \ll_\epsilon (\lambda+\epsilon)^k.
\end{equation}
From numerical data, $\lambda \approx 0.78$.
\end{conj}

The exponential bounds of Theorem~\ref{thm: bounds} already imply that the $r$\textsuperscript{th} moment remains bounded for any $r\geq 1$.
\begin{cor}\label{cor: moments}
The limit of the $r$\textsuperscript{{\rm th}} moment of $M_{[0,n-1]}$,
\begin{equation}
\lim_{n\rightarrow \infty} \Expect{M_{[0,n-1]}^r},
\end{equation}
exists and is finite.
\end{cor}

\begin{thm}\label{thm: variance}
The limit
\begin{equation}
\lim_{n\rightarrow \infty} \Var \left(M_{[0,n-1]}\right)
\end{equation}
exists and is about $35.9658$, as these are the first digits of its decimal expansion. This limit can be written as the following
convergent series with exponential decay:
\begin{equation}\label{eqn: varconvergent}
\lim_{n\rightarrow \infty} \Var \left(M_{[0,n-1]}\right)\ =\ 4\lim_{n\rightarrow \infty} \sum_{i < j < n} \pp(i \mbox{ and } j \not \in A+A ) - 40.
\end{equation}
%see Proposition \ref{prop: exact} for the exact (and long) formula for $\pp(i \mbox{ and } j \not \in A+A )$.
\end{thm}

Note that ``$i \mbox{ and } j \not \in A+A$'' is meant to be parsed as ``$(i\not\in A+A) \text{ AND } (j\not\in A+A)$''.

%%%%%%%%%%%%%%%%%%%%%%%%%%%%%%%%%%%%%%%%%%%%%%%%%%%%%%%%%%%%%%%%%%%%%%%%%%%%%%%%%%%%%%%%%%%%%%%%%%%%%%%%%
\subsection{Variance and Decay Rates of Missing Sums}\label{sec:subsecvariancedecayrates}

The bounds in Theorem \ref{thm: bounds} are due to formulas for probabilities of events such as
\begin{equation}
\pp(a_1,\ a_2,\ \dots,\ \mbox{and}\ a_m  \not \in A+A),\end{equation}
by which we mean the probability that all of $a_1,a_2,\dots,a_m$ are in the complement of $A+A$.
This represents the probability that a particular configuration is not in $A+A$. As long as $n > a_m$, there is no dependence on $n$ since this probability just depends on $[0, a_m] \cap A$. We therefore can assume that $A \subseteq [0,a_m]$.
Formulas for such probabilities are also important for finding the moments of $M_{[0,n-1]}$. For example, to find the expectation of $|A+A|$, \cite{MO} find an exact formula for $\pp(k \not \in A+A)$, which is approximately
\begin{equation}\label{eqn: iapp}
\pp(k \not \in A+A)\ =\ \Theta((3/4)^{k/2}),
\end{equation}
where we say $g(n) = \Theta(f(n))$ if there exist constants $C_1,C_2$ such that for all $n$
\begin{equation}
C_1 f(n)\ \le\ g(n)\ \le\ C_2 f(n).
\end{equation}
Similarly, to find the variance, we can study $\pp(i \mbox{ and } j \not \in A+A)$ as seen from the series expansion  in
\eqref{eqn: varconvergent}.
In Proposition \ref{prop: exact}, we find an exact formula for this probability and in Corollary \ref{cor: ij}, we show that for fixed $m$ we have the following approximation:
\begin{equation}\label{eqn: ijapp}
\pp(k \mbox{ and } k+ m  \not \in A+A)\ =\ \Theta((\phi/2)^k).
\end{equation}
The implied constants in \eqref{eqn: ijapp} depend significantly on $m$ and in Corollary \ref{cor: ij}, we also find these constants.

Note that both \eqref{eqn: iapp} and \eqref{eqn: ijapp}  are exponential in $k$. In fact, we prove that in general such probabilities are approximately exponential in $k$.

\begin{thm}\label{thm: configuration}
For any fixed $a_1, \dots, a_m$, there exists $\lambda_{a_1, \dots, a_m}$ such that
\begin{equation}
\pp(k + a_1, k + a_2, \dots, \mbox{ and } k + a_m \not \in A+A)\ =\ \Theta(\lambda_{a_1, \dots, a_m}^k),
\end{equation}
where the implied constants depend on $a_1, \dots, a_m$ but not $k$.
\end{thm}

The fact that $\pp(k + a_1, k + a_2, \dots, \mbox{ and } k + a_m \not \in A+A)$ is approximately exponential
supports Conjecture \ref{conj: main} that the distribution of missing sums is approximately exponential.

For the particular configuration $a_1 = 1, a_2 = 2, \dots, a_m = m$, the case of consecutive missing elements, we can approximate
$\lambda_{a_1, \dots, a_m}$ well as seen in the following theorem.

\begin{thm}\label{thm: consecutive}
For any $k,m$
\begin{equation}
\left(\frac{1}{2}\right)^{(k+m)/2}
\ \ll\
\pp(k + 1, k +2, \dots, \mbox{ and } k +m \not \in A+A)
\ \ll\ \left(\frac{1}{2} \right)^{(k+m)/2}(1+ \epsilon_m)^k,
\end{equation}
with $\epsilon_m \rightarrow 0$ as $m \rightarrow \infty$.
To be more precise, the exact form of upper bound is
$(1/2)^{(k+m)/2} 2^{k/m} $. This implies that
\begin{equation}
\lambda_{0, 1, 2, \dots, m} \rightarrow \left(\frac{1}{2}\right)^{1/2}
\end{equation}
as $m \rightarrow \infty$.
\end{thm}

As we will see in the proof of Theorem \ref{thm: bounds}, the lower bound $(1/2)^{(k+m)/2}$ is essentially the probability of missing the first $k+m$ elements in $A+A$. By Theorem \ref{thm: consecutive}, we have that for large $m$, $\pp(k + 1, k +2, \dots, \mbox{ and } k +m \not \in A+A)$ is also approximately $(1/2)^{(k+m)/2}$. This means that for large $m$, essentially the only way to miss $m$ consecutive elements in $A+A$ starting at $k+1$ is through the trivial way - namely missing all of the first $k+m$ elements of $A+A$.

Theorem \ref{thm: consecutive} is in fact a special case of the following inequality.
\begin{thm}\label{thm: lambdabounds}
For $\lambda_{a_1, \dots, a_m}$ with $0 \le a_1 < \cdots < a_m$,
\begin{equation}
\lambda_{a_1, \dots, a_m} \le \pp(A,B \subseteq [0, \lfloor a_m/2 \rfloor ] \mid a_1, \dots, a_m \not \in A+B)^{\frac{1}{a_m +2}}.
\end{equation}
where $A, B$ are two independently chosen sets.
\end{thm}

%%%%%%%%%%%%%%%%%%%%%%%%%%%%%%%%%%%%%%%%%%%%%%%%%%%%%%%%%%%%%%%%%%%%%%%%%%%%%%%%%%
\subsection{Other types of random sets and the divot}\label{subsec: intro divot}

Figure~\ref{figure:x_k Experimental graph} shows a surprising phenomenon: experimentally, $$m(7)<m(6)<m(8).$$ That is, a random subset of $[0,10^{10}]$ is more likely to have a sumset missing 6 (or 8) elements than one missing 7 elements. That is, the distribution of $M_{[0, n-1]}$ appears to be bimodal for large $n$. We have made a massive computation (details in \S\ref{sec:bimodality}), looping over $2^{43}$ sets and using only 64-bit integer arithmetic, that lead to the following bounds:
\be 0.07177< m(6) < 0.07202,\ \ 0.07138 < m(7) < 0.7170,\ \ 0.07243 < m(8) < 0.07282.\ \ \ \ee
We note that our bounds are actually in the form
\begin{multline*}
 \frac{107418021089142422011644549535908507304608994344051}{1496577676626844588240573268701473812127674924007424}
 < m(6) \\ m(6) <
 \frac{620778536995376440633741122321102716502820362028980739}{8620287417370624828265702027720489157855407562282762240};
\end{multline*}
we hope the reader will excuse our preference for reporting equivalent decimals, rounded in the proper directions to maintain truth.

Closer inspection of Figure~\ref{figure:x_k Experimental graph} also reveals an apparent parity effect: $$m({2k})+m({2k+2}) > 2 m({2k+1}).$$ Here are two plausible explanations for this. The first is that $M_{[0,n-1]}$ is essentially the sum of two iidrvs: the number of missing sums in $[0,n-1]$ and in $[n,2n-2]$. For any two iidrvs $X_1,X_2$ taking integer values, $\Prob{X_1+X_2\text{ even}}\geq \Prob{X_1+X_2\text{ odd}}$, as the calculation comes down to $x^2+y^2 \geq 2xy$. Another parity effect is observed on the ends: as soon as $0\not\in A$, then both $0$ and $1$ are not in $A+A$. Thus, on the ends, $A+A$ always misses an even number of sums.

To compensate for these observations, it is necessary to consider the connections between different ways of selecting a random set. We consider uniformly selecting subsets of $[0,n-1]$, subsets of $[0,n]$ with diameter $n$, subsets of $\NN$, and subsets of $\NN$ that contain 0. We lay out our notation as follows:\\

\begin{center}
\begin{tabular}{|ccccc|} \hline
   set & setting & condition & missing sums & $\pp$(missing $k$ sums) \\ \hline
   $A$ & $[0,n-1]$ & $\emptyset$ & $M_{[0,n-1]} := 2n-1-|A+A|$ & $m_n(k)$ \\
   $B$ & $[0,n]$ & $\{0,n\}\subseteq B$ & $M_{[0,n]\mid \{0,n\}} := 2n+1-|B+B|$ & $w_n(k)$\\
   $C$ & \NN & $\emptyset$ & $M_{\NN} := | \NN \setminus (C+C)|$ & $y(k)$ \\
   $D$ & \NN & $0\in D$ & $M_{\NN \mid \{0\}}:= |\NN \setminus (D+D)|$ & $z(k)$ \\ \hline
\end{tabular}
\end{center}
\noindent Additionally, we set $m(k):=\lim_{n\to\infty} m_n(k)$ and $w(k):= \lim_{n\to\infty} w_n(k)$.\\ \

\begin{figure}
\begin{center}
\includegraphics[scale=0.6]{mkGraph.eps} \includegraphics[scale=0.6]{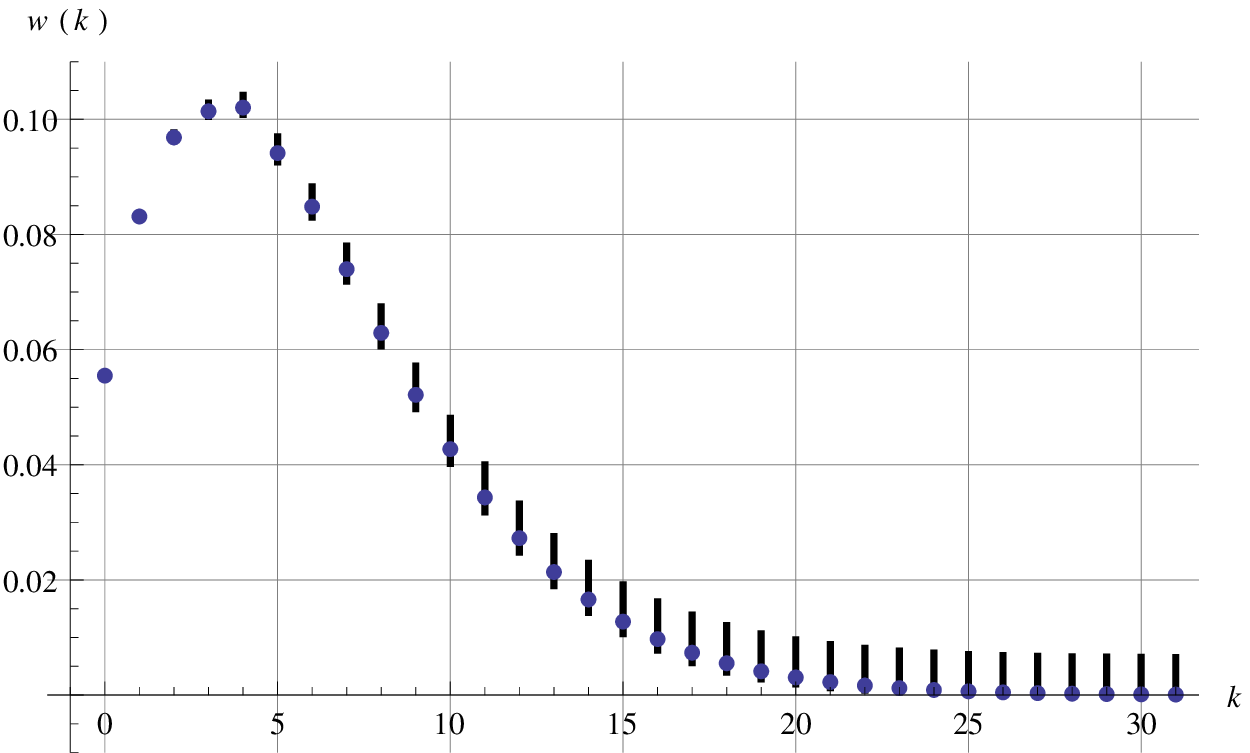}

\includegraphics[scale=0.6]{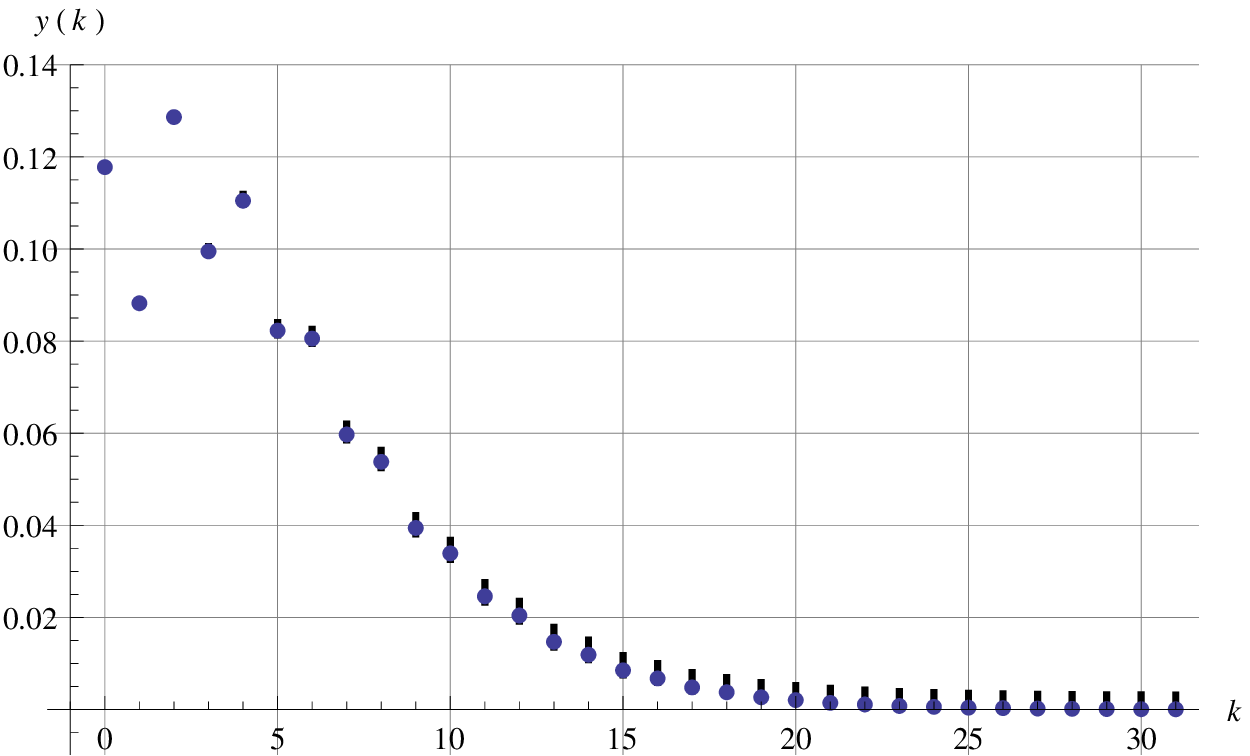} \includegraphics[scale=0.6]{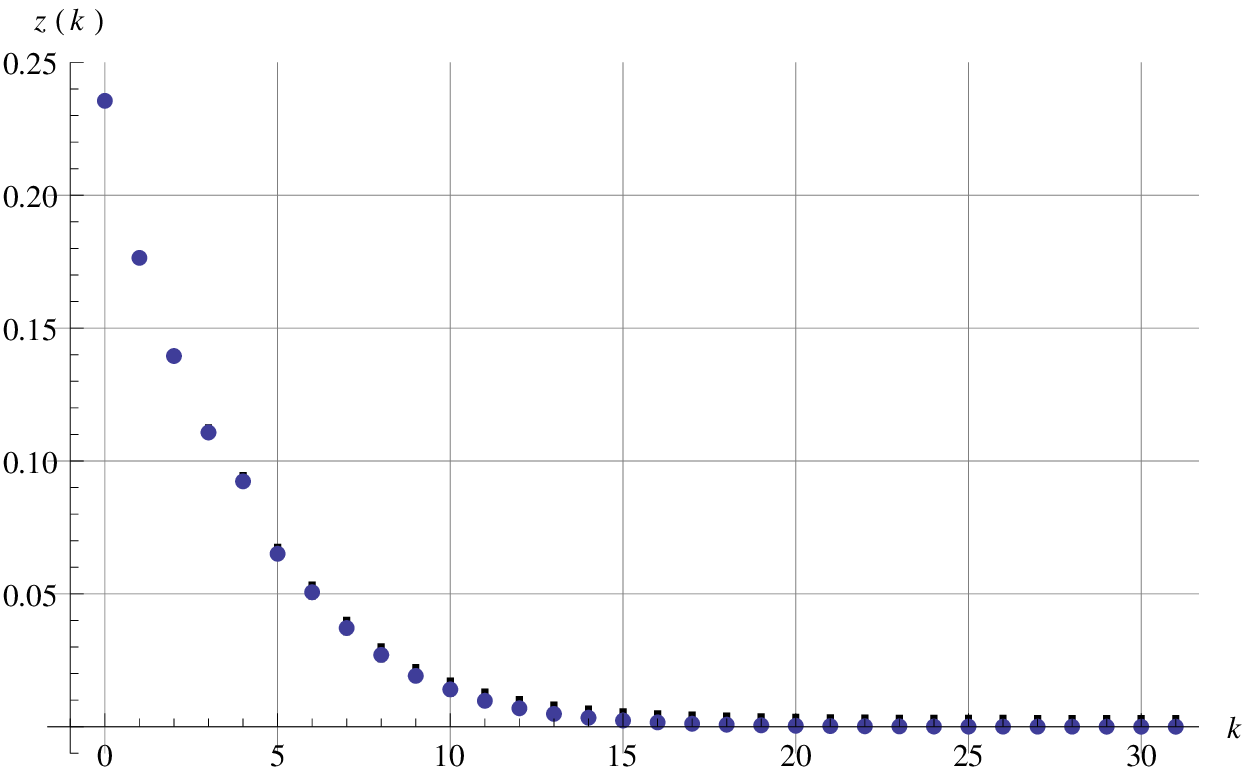}
\end{center}
\caption{Experimental values of $m(k)$, $w(k)$, $y(k)$, $z(k)$, with vertical bars depicting the values allowed by our rigorous bounds. See \S\ref{sec:bimodality} for details.\label{figure:z_k Experimental graph}}
\end{figure}

Our first parity-effect observation essentially boils down to
   \begin{equation}\label{equ: m and y}
      m_n(k) \to \sum_{i=0}^k y(i)y(k-i),
   \end{equation}
a rigorous exposition of this can be found in~\cite{In} and is sketched in \S\ref{subsec: connecting the distribs}. The second observation and Bayes' Theorem leads us to
   \begin{equation} \label{equ: y and z}
      y(k) = \sum_{i=0}^\infty \Prob{\min C = i}\Prob{|[2i,\infty) \setminus (C+C)| = k-2i} = \sum_{i=0}^{\lfloor k/2\rfloor} 2^{-(i+1)} z(k-2i).
   \end{equation}
Similarly to~\eqref{equ: m and y}, one can prove that
   \begin{equation} \label{w and z}
      w_n(k) \to \sum_{i=0}^k z(i) z(k-i).
   \end{equation}
Thus, all four distributions can be understood in terms of $z(k)$. Experiments and our bounds (see Figure~\ref{figure:z_k Experimental graph} for small values of $k$) indicate that $M_{\NN \mid \{0\}}$ has an approximately geometric distribution, and exhibits no obvious parity effect. Computationally, we focus on bounding $z$ and then allow this to determine bounds on $m$, $w$ and $y$.

We bound $z$ by conditioning on $I:=D\cap[0,44)$, and loop over all $2^{43}$ possible values of $I$ (a priori, $0\in I$). For each $I\subseteq[0,44)$, we explicitly know $(D+D)\cap[0,44)$, we have much information concerning $(D+D)\cap[44,88)$, and theoretically $(D+D)\cap[88,\infty)$ is $[88,\infty)$ with high probability. This allows us to give reasonable upper and lower bounds on $\Prob{M_{\NN \mid \{0\}}=k \mid D\cap[0,44)=I}$ for each $I$.

If we suppose that $M_{\NN \mid \{0\}}$ is exactly geometric with parameter $\lambda$ (i.e., set $z(k)=(1-\lambda)\lambda^k$) and define $y(k)$ and $m(k)$ using \eqref{equ: m and y} and \eqref{equ: y and z}, we find that the distribution of $M_{\NN \mid \{0\}}$ would be bimodal with a divot at $k=7$ only for the narrow parameter range $0.756<\lambda<0.771$. The best-squares fit for $\lambda$ is $0.765$.
%\marginpar{\tiny This was fitting the values, but it would be more appropriate to fit the logarithms.}
If we suppose that $M_{\NN \mid \{0\}}$ has a Poisson distribution, i.e., $z(k)=\lambda^k e^{-\lambda}/k!$, we find that there are no $\lambda$ whatsoever that give a bimodal distribution with divot at $k=7$.

This implies that the divot's existence relies not only on the above observations but also on the specific values of $z_k$ for small values. We note that $z_4$ in particular is larger than the geometric model predicts; more than half of the least-squares error is from $z_4$. The rigorous bounds we give also show this bias towards 4, though we currently have no understanding as to why this is the situation.

\begin{thm}
The limits defining $m(k)$ and $w(k)$ are well-defined, positive, and $\sum_{k=0}^\infty m(k) =\sum_{k=0}^\infty w(k) =1$. Rigorous bounds on $m(k)$, $w(k)$, $y(k)$ and $z(k)$ for $0\leq k < 32$ are given in Appendix~\ref{Appendix}. In particular, $m(7)<m(6)<m(8)$.
\end{thm}

%%%%%%%%%%%%%%%%%%%%%%%%%%%%%%%%%%%%%%%%%%%%%%%%%%%%%%%%%%%%%%%%%%%%%%%%%%%%%%%%%%
%%%%%%%%%%%%%%%%%%%%%%%%%%%%%%%%%%%%%%%%%%%%%%%%%%%%%%%%%%%%%%%%%%%%%%%%%%%%%%%%%%
\section{Graph-Theoretic Framework}\label{sec: graph}

We first develop a graph-theoretic framework to study dependent random variables and calculate probabilities like
$\pp(a_1, \dots, \mbox{ and } a_m \not \in A+A)$. Note that for odd $i$
\begin{equation}
\{i\not \in A+A\}\ =\
\left\{ \left( 0\not \in A \mbox{ or } i\not \in A\right) \mbox{ and } \cdots \mbox{ and } \left( (i-1)/2\not \in A \mbox{ or }(i+1)/2 \not \in A \right) \right\},
\end{equation}
and for even $i$
\begin{equation}
\{i\not \in A+A\}\ =\ \{ \left( 0\not \in A \mbox{ or } i\not \in A\right) \mbox{ and }\cdots \mbox{ and } \left( i/2-1\not \in A \mbox{ or } i/2+1 \not \in A\right)  \mbox{ and } i/2 \not \in A\}.
\end{equation}
Therefore for distinct $i$ the events $\{i \not \in A+A\}$ are dependent as both depend on conditions on $A$ like
$\{ 0 \not \in A\} $.

For example, the conditions on $A$ necessary for $\{3 \mbox{ and } 7 \not \in A+A\}$ are
\begin{equation}\label{eqn: conditionarray}
\begin{array}{lllll}
& i = 3: & 0 \mbox{ or } 3 \not \in A & \: \: \: \: \: \: \: j = 7:  & 0 \mbox{ or } 7 \not \in A\\
& &\mbox{and } 1 \mbox{ or } 2 \not \in A  & & \mbox{and } 1 \mbox{ or } 6 \not \in A\\
& & & & \mbox{and } 2 \mbox{ or } 5 \not \in A\\
& & & & \mbox{and } 3 \mbox{ or } 4 \not \in A.
\end{array}
\end{equation}
Since the two lists have integers in common, there is dependence between the events $\{3\not \in A+A\}$ and $\{7\not \in A+A\}$.

We construct a graph to represent the dependencies between the random variables. We call this graph the
\emph{condition graph} for the probability. We construct the condition graph for  $\pp(a_1, \dots$, and $a_m$ $\not$ $\in$ $A+A)$, where
$a_1 < \cdots < a_m$, in the following way:
\begin{enumerate}
\item For every integer in $[0, a_m]$, add a vertex labeled with that integer.
\item Add an edge between two vertices labeled with $i$ and $j$ if $i + j = a_k$ for some $ 1\le k \le m$.
\end{enumerate}
See Figure \ref{fig: condition37} for the condition graph for $\pp(3 \mbox{ and } 7 \not \in A+A)$.

\begin{figure}[h]
\begin{tikzpicture}
  [scale=.7,auto=left,every node/.style={circle,fill=black!20,minimum size=12pt,inner sep=2pt}]
  \node (n0) at (0,2)       {0};
  \node (n1) at (1.4,1.4)   {1};
  \node (n2) at (2,0)       {2};
  \node (n3) at (1.4,-1.4)  {3};
  \node (n4) at (0,-2)      {4};
  \node (n5) at (-1.4,-1.4) {5};
  \node (n6) at (-2,0)      {6};
  \node (n7) at (-1.4,1.4)  {7};

  \foreach \from/\to in {n0/n7,n1/n6,n2/n5,n3/n4,n0/n3,n1/n2}
    \draw (\from) -- (\to);
\end{tikzpicture}
\caption{Condition Graph for $\pp(3 \mbox{ and } 7 \not \in A+A)$.}
\label{fig: condition37}
\end{figure}

By construction, we have a one-to-one correspondence between edges and conditions and vertices and integers in $[0, a_m]$.
For example, the edge between vertices labeled with $1$ and $6$ represents the condition that $1$ or $6 \not \in A$, which is one of the conditions necessary for $7 \not \in A+A$ in \eqref{eqn: conditionarray}. For each condition, we need to pick at least one element to exclude from $A$. Therefore in the condition graph, for each edge we need to pick at least one of its vertices.
That is, we need to pick a vertex cover (recall a vertex cover of a graph is a set of vertices such that each edge is incident to at least one vertex in the set). Using this method, we get the following lemma.

\begin{lem}\label{lem: vertexcover}
$\pp(a_1, \dots, \mbox{ and } a_m \not \in A+A)$ equals the probability that  we chose a vertex cover for the condition graph.
\end{lem}

Note that when we pick vertices in the condition graph for our vertex cover, we are picking to exclude those vertices from $A$.
For example, note that the vertices $7, 0, 4$ and $6, 2$ form a vertex cover for the condition graph of $\pp(3 \mbox{ and } 7 \not \in A+A)$ in Figure \ref{fig: condition37}. Then if $7,0,4,6,2 \not \in A$, then $3 \mbox{ and }7 \not \in A+A$ since all conditions in \eqref{eqn: conditionarray} are met.

Finally, note that when we calculate the probability of chosing a vertex cover for the condition graph,
we no longer need to consider a labeled graph. This is because vertices represent elements of $A$, and since each element of $A$ is equally likely to be chosen
(as $A$ is chosen uniformly randomly), we do not need to differentiate between different elements.

%%%%%%%%%%%%%%%%%%%%%%%%%%%%%%%%%%%%%%%%%%%%%%%%%%%%%%%%%%%%%%%%%%%%%%%%%%%%%%%%%%%%%%%%%%%%%%
%%%%%%%%%%%%%%%%%%%%%%%%%%%%%%%%%%%%%%%%%%%%%%%%%%%%%%%%%%%%%%%%%%%%%%%%%%%%%%%%%%%%%%%%%%%%%%
\section{Variance of Missing Sums}\label{sec: variance}

We now use the graph-theoretic framework from the previous section to prove Theorem \ref{thm: variance} and find the variance.
%We will abbreviate $M_{[0,n-1]}$ to $M_{[0, n-1]}$ henceforth.

We first note that the result of \cite{MO} in \eqref{thm: martinobryant} is really that
\begin{equation}\label{eqn: martin2}
\E\left[M_{[0, n-1]}(A)\right]\ =\  \sum_{0\le i\le 2n-2} \pp( i \not\in A+A)  = 10 + O((3/4)^{n/2}).
\end{equation}
Since
\begin{equation}\label{eqn: varform}
\Var\left(M_{[0, n-1]}(A)\right)\ =\ \E\left[M_{[0, n-1]}(A)^2\right] - \left(\E\left[M_{[0, n-1]}(A)\right]\right)^2
\end{equation}
and we know $\E\left[M_{[0, n-1]}(A) \right]$ from \eqref{eqn: martin2},
to find the variance we just need to determine $\E\left[M_{[0, n-1]}(A)^2\right]$, which equals the following:
\begin{eqnarray}\label{eqn: var}
\E\left[M_{[0, n-1]}(A)^2\right] & \ = \ & \frac{1}{2^n} \sum_{A\subseteq [0,n-1]}  |\{\mbox{missing sums in } A+A \}|^2 \nonumber\\
&=& \frac{1}{2^n} \sum_{A\subseteq [0,n-1]} \sum_
{\substack{
0\le i,j \le 2n-2\\
i,j \not \in A+A
}} 1 \nonumber\\
&=& \frac{1}{2^n} \sum_{0\le i,j\le 2n-2} \sum_
{\substack{
A\subseteq [0,n-1]\\
 i,j \not \in A+A
}} 1\nonumber\\
&=& \sum_{0\le i,j\le 2n-2} \pp(A \subseteq [0,n-1] \mid i \mbox{ and }j \not\in A+A)\nonumber\\
&=& 2\sum_{0\le i < j\le 2n-2} \pp( i \mbox{ and }j \not\in A+A) +
 \sum_{0\le i\le 2n-2} \pp( i \not\in A+A).
\end{eqnarray}
Combining \eqref{eqn: varform}, \eqref{eqn: martin2}, and \eqref{eqn: var}, we get
\begin{equation}\label{eqn: variancenotsimple}
\Var\left(M_{[0, n-1]}(A)\right)\ =\ 2\sum_{0\le i < j\le 2n-2} \pp( i \mbox{ and }j \not\in A+A) - 90 +  O((3/4)^{n/2}).
\end{equation}

We first simplify the sum over $i,j$. Note that if $i, j < n$, then
\begin{equation} \label{eqn: flip}
\pp(i \mbox{ and }j \not \in A+A)\ = \ \pp( 2n-2- i \mbox{ and } 2n-2- j \not \in A+A),
\end{equation}
and so
\begin{equation}\label{eqn: flipsum}
\sum_{0 \le i < j < n}\pp(i \mbox{ and }j \not \in A+A)\ =\ \sum_{n \le i < j \le 2n-2}\pp(i \mbox{ and }j \not \in A+A).
\end{equation}
Also, note that if $i < n/2$ and $j > 3n/2$, then $\{ i \not \in A+A\}$ and $\{ j \not \in A+A\}$ are independent. This is because $\{ i \not \in A+A\}$ depends only on $[0,i]\cap A$ and $\{ j \not \in A+A\}$ depends only on $[j-n+1, n-1] \cap A$ and if $i < n/2$ and $j > 3n/2$, these sets are disjoint.
Therefore for such $i,j$, we have
\begin{equation}
\pp(i \mbox{ and }j \not \in A+A)\ =\ \pp(i \not \in A+A) \pp(j \not \in A+A).
\end{equation}
Finally note that if $n/2 \le i < n$ or  $n \le j \le 3n/2$, then
\begin{equation}
\pp( i\mbox{ and }j \not \in A+A)\ =\ O((3/4)^{n/4})
\end{equation}
by \eqref{eqn: iapp}.
Therefore
\begin{eqnarray}\label{eqn: opposite}
&& \sum_{i < n,\  n \le j} \pp( i\mbox{ and }j \not \in A+A) \nonumber\\
&&=  \sum_{i < n/2 {\rm\ and\ }  3n/2 < j} \pp(i\mbox{ and }j \not \in A+A)
+ \sum_{ n/2 \le i < n {\rm\ or\ }  n \le j \le 3n/2} \pp(i\mbox{ and }j \not \in A+A)  \nonumber\\
&&=  \sum_{i < n/2,\  3n/2 < j} \pp(i\mbox{ and }j \not \in A+A) + O(n^2 (3/4)^{n/4})\nonumber\\
&&= \left( \sum_{i < n/2} \pp( i \not \in A+A) \right) \cdot \left( \sum_{3n/2 < j \le 2n-2}\pp( j \not \in A+A)\right)   + O(n^2 (3/4)^{n/4})\nonumber\\
&&=  \left(5+O((3/4)^{n/4})\right) \cdot \left(5+O((3/4)^{n/4})\right)  + O(n^2 (3/4)^{n/4})\nonumber\\
&&= 25 + O(n^2(3/4)^{n/4}),
\end{eqnarray}
where we  use \eqref{eqn: martin2} and \eqref{eqn: flip}  to get the second to last equality.
Combining \eqref{eqn: flipsum} and \eqref{eqn: opposite}, we have
\begin{eqnarray}
&& \sum_{0 \le i < j \le 2n-2 }\pp( i \mbox{ and } j \not \in A+A) \nonumber\\
&&\:\: =\sum_{0 \le i < j < n}\pp(i \mbox{ and } j \not \in A+A)
+\sum_{n \le i < j \le 2n-2}\pp(i \mbox{ and } j \not \in A+A)
+ \sum_{ i <  n, \ n \le j}\pp(i \mbox{ and } j \not \in A+A)\nonumber\\
&&\: \: = 2\sum_{0\le i < j \le n-1}\pp(i \mbox{ and } j \not \in A+A) + 25 + O(n^2(3/4)^{n/4}),
\end{eqnarray}
and so by \eqref{eqn: variancenotsimple}
\begin{equation}\label{eqn: variancesimple}
\Var\left(M_{[0, n-1]}(A)\right)\ =\ 4\sum_{0\le i < j \le n-1}\pp(i \mbox{ and } j \not \in A+A) -40+ O(n^2(3/4)^{n/4}).
\end{equation}
Therefore to find the variance, we just need to study $\pp(i \mbox{ and }j \not \in A+A)$ for $i <  j < n$.

\emph{Since the other cases are handled similarly, we only present the details for the case when $i$ and $j$ are both odd.}
By Lemma \ref{lem: vertexcover}, we just need to study the condition graph for $\pp(i\mbox{ and }j \not \in A+A)$. Recall that we already found the condition graph for $\pp(3 \mbox{ and }7 \not \in A+A)$ in Figure \ref{fig: condition37}. After untangling this graph, we see that it really consists of two components, as seen in Figure \ref{fig: untangled37}.

\begin{figure}[h]
\begin{tikzpicture}
[scale=.7,auto=left,every node/.style={circle,fill=black!20,minimum size=12pt,inner sep=2pt}]
  \node (n7) at (4,1)  {7};
  \node (n0) at (6,1)  {0};
  \node (n3) at (8,1)  {3};
  \node (n4) at (10,1) {4};
  \node (n6) at (14,1)  {6};
  \node (n1) at (16,1)  {1};
  \node (n2) at (18,1)  {2};
  \node (n5) at (20,1)  {5};

  \foreach \from/\to in {n5/n2,n2/n1,n1/n6,n4/n3,n3/n0,n0/n7}
  \draw (\from) -- (\to);
\end{tikzpicture}
\caption{Untangled condition graph for $\pp(3 \mbox{ and }7 \not \in A+A)$.}
\label{fig: untangled37}
\end{figure}

Also note that each component is a \emph{segment graph}, a graph that consists of a sequence of vertices such that each vertex is connected only to the vertices to its immediate left and right.  A similar situation holds in general, as seen by the following proposition.

\begin{prop}
The condition graph for $\pp( i \mbox{ and } j \not \in A+A)$ has components that are segment graphs.
\end{prop}

\begin{proof}
The condition graph for $\pp( i \mbox{ and } j \not \in A+A)$ has vertices with degree less than or equal to $2$; if the vertex is labeled with $\ell$, it can only be connected to vertices labeled $i - \ell$ or $j-\ell$ (if such vertices exist).

Furthermore, there are no cycles in the condition graph. Suppose there is a cycle in the condition graph. Consider the vertex in the cycle with the maximum label $\ell$ and consider the vertices around this vertex. Each of these vertices must have exactly two neighbors and so we have the following situation as seen in Figure \ref{fig:figuresneedlabels}.
\begin{figure}[h]\tiny
\begin{tikzpicture}
 [scale=1.,auto=left,every node/.style={circle,fill=black!20,minimum size=22pt,inner sep=1pt}]
  \node (n1) at (4,1)  {\:\:$\ j-\ell\ $\:\:};
  \node (n2) at (7,1)  {\:\:\:\:\:$\ \ \ell\  \ $\:\:\: \:};
  \node (n3) at (10,1)  {\:\:$\ \ i-\ell\ \ $\:\:};
  \node (n4) at (13,1)  { $\ \ell + j - i\ $};
  \foreach \from/\to in {n1/n2,n2/n3,n3/n4}
    \draw (\from) -- (\to);
\end{tikzpicture}
\caption{Vertices around a labeled vertex $\ell$.}
\label{fig:figuresneedlabels}
\end{figure}
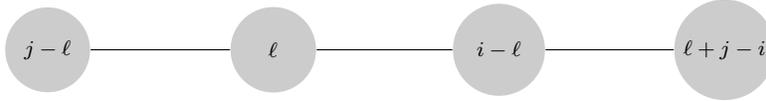 \normalsize

Notice that $\ell+j-i > \ell $ since $j> i$. Therefore, $\ell$ is not the maximum label, which is a contradiction and proves that we cannot have a cycle. Thus all components are trees with all vertices having degree less than or equal to $2$, implying that all are segment graphs.
\end{proof}

Since labels in different components are distinct and there are no edges between different components, each component is independent. That is, the probability of getting a vertex cover for the entire graph is the product of the probability of getting vertex covers for each component. In this way, we just need to find the probability of getting a vertex cover for each component. To do this, we find the number of vertex covers for an arbitrary segment graph, which we do in the following proposition.

\begin{prop}
The number of vertex covers $g(n)$ for a segment graph with $n$ vertices satisfies
$g(n) = F_{n+2}$, where $F_k$ is the $k$\textsuperscript{th} Fibonacci number.
\end{prop}

\begin{proof}
There are two cases: the first vertex of the segment graph is in the vertex cover, or it is not. If the first vertex is in the cover, then the first edge already has one of its vertices picked. Therefore we just need a vertex cover for the subgraph with $n-1$ vertices that follows the first edge, and by definition there are $g(n-1)$ such covers. If the first vertex is not in the cover, then the second vertex must be the cover since the first edge must have one of its vertices chosen. Since the second vertex is now in the cover, then the second edge automatically has one of its vertices in the cover. Therefore we just need a vertex cover for the subgraph with $n-2$ vertices that follows the second edge, and by definition there are $g(n-2)$ such vertex covers. Therefore, we have the Fibonacci recursive relationship $g(n)= g(n-1)+ g(n-2)$. As $g(2)=3 = F_4$ and $g(3)= 5 = F_5$, these initial conditions and the recurrence imply $g(n) = F_{n+2}$, completing the proof.
\end{proof}

Therefore, we have
\begin{equation}
\pp(\mbox{chose a vertex cover for a segment graph with $n$ vertices})\ =\ \frac{F_{n+2}}{2^n}.
\end{equation}

Returning to our example with $3$ and $7$, we note that since the condition graph in this case consists of two segment graph components each of length $4$, we have
\begin{equation}
\pp(3 \mbox{ and } 7 \not \in A+A)\ =\ \frac{F_{4+2}}{2^4} \cdot \frac{F_{4+2}}{2^4}\ =\ \frac{1}{4},
\end{equation}
where we can multiply the probabilities by the independence of the components.

In general, as the condition graph may have many components we must find how many segment graph components there are in the entire graph for $\pp(i \mbox{ and }j\not \in A+A)$.

\begin{prop}\label{prop: compnum}
There are $(j-i)/2$ segment graph components for the graph of $\pp(i \mbox{ and } j \not \in A+A)$.
\end{prop}

\begin{proof}
Note that in total $j+1$ vertices are used in the graph of $\pp(i \mbox{ and } j \not \in A+A)$; since
$\{i \mbox{ and } j \not \in A+A\}$ depends just on $A\cap [0, j]$, the graph uses exactly the integers in $[0,j]$. Also note that each component must end with a vertex labeled by an integer greater than $i$. If a component ends with a vertex labeled by
$\ell \le i$, then it can be connected to two other vertices $i - \ell$ and $j - \ell$. Remember that we are assuming $i$ and $j$ are odd (the other cases are similar). As they are odd, $i - \ell  \ne \ell$ and $j - \ell \ne \ell$ and so $i-\ell, j - \ell, \ell$ are all distinct. Since $\ell$ is connected to two other vertices, it cannot be an end vertex. Therefore, each end vertex is labeled by some integer in $[i+1,j]$. Also note that each of these integers must be end vertex since it cannot be used to add up to $i$. Therefore, the set $[i+1, j]$ coincides with the set of end vertices and since each component has two end vertices with distinct labels, there are $(j - i)/2$ components.
\end{proof}

We also need to find the length of each component. Fortunately, there are only two possible component lengths for the graph of
$\pp(i \mbox{ and } j \not \in A+A)$, as seen by the following lemma.

\begin{prop}\label{prop: complength}
The length of each segment graph component for the graph of $\pp(i \mbox{ and } j \not \in A+A)$ is always either
\begin{equation}\label{eq:twovaluesfromprop}
2\left\lceil \frac{i+1}{j-i} \right\rceil \: \: \: {\rm or} \: \: \: 2\left\lceil \frac{i+1}{j-i} \right\rceil+2.
\end{equation}
\end{prop}

\begin{proof}
First note that the difference between a given vertex and another vertex that is two edges away is
$j-i$. This is because the sum of the vertices that share an edge alternates between $i$ and $j$, so that we have segments of the form given in Figure \ref{fig:figiminusxxjminusx}. The difference between $j-x$ and $i-x$ is $j-i$ as needed.

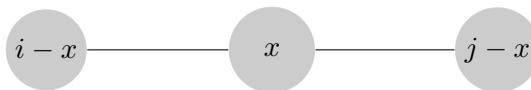
\begin{figure}[h]\small
\begin{tikzpicture}
 [scale=1.,auto=left,every node/.style={circle,fill=black!20,minimum size=20pt,inner sep=0.5pt}]
  \node (n1) at (4,1)  {\:$i-x$\:};
  \node (n2) at (7,1)  {\:\:\:\: $x$\ \:\:\:\:};
  \node (n3) at (10,1)  {\:$j-x$\:};
  \foreach \from/\to in {n1/n2,n2/n3}
    \draw (\from) -- (\to);
\end{tikzpicture}
\caption{Difference between every other vertex.}
\label{fig:figiminusxxjminusx}
\end{figure}\normalsize

Now note that these differences can be used to determine the size of each component. Suppose the end vertex of a segment graph is $m$. Since we decrease by $j-i$ for every two vertices and since we only use non-negative integers, there can only be
\begin{equation}
\left \lfloor \frac{m}{j-i}\right \rfloor + 1\ =\ \left \lceil \frac{m+1}{j-i} \right \rceil
\end{equation}
decreases. Since we decrease once for every two vertices, we have that the length is twice the number of decreases. Therefore the length is
\begin{equation}
 2\left \lceil \frac{m+1}{j-i} \right \rceil.
\end{equation}

From Proposition \ref{prop: compnum}, we also know that the end vertex $m$ of each segment graph satisfies $i < m \le j$. Therefore, the length of each segment graph is always
\begin{equation}\label{eqn: k}
2\left \lceil \frac{i+1}{j-i} \right \rceil\ \mbox{ or }\
 2\left \lceil \frac{i+1}{j-i} \right \rceil +2,
\end{equation}
as desired.
\end{proof}

For simplicity, we denote the first of the two values in \eqref{eq:twovaluesfromprop} by $q$ and the second by $q+2$. We must find the number of components with size $q$ and $q+2$. Suppose there are $r$ components of size $q$ and $r'$ components of size $q+2$. Then conditions on the number of components from Proposition \ref{prop: compnum} and the length of each component from Proposition \ref{prop: complength} gives us the following two equations:
\begin{eqnarray}\label{eqn: system}
qr + (q + 2)r' &\ = \ & j+1\nonumber\\
r+r' &=& \frac{j-i}{2}.
\end{eqnarray}
Solving these equations for $r, r'$ in terms of $q$ gives
\begin{eqnarray}\label{eqn: m}
r  &\ = \ & \frac{1}{2}\left(\frac{j-i}{2}q - (i +1)\right)
\ = \ \frac{1}{2}\left((j-i)\left\lceil \frac{i+1}{j-i} \right\rceil - (i +1)\right) \nonumber\\
r' &\ = \ & \frac{1}{2}\left(j+1 -\frac{j-i}{2}q \right)\  \ \ \ =\ \frac{1}{2}\left(j+1 -(j-i)\left\lceil \frac{i+1}{j-i} \right\rceil \right).
\end{eqnarray}
Therefore, again by independence of components, we have for odd $i,j$ that
\begin{eqnarray}
\pp(i \mbox{ and } j \not \in A+A)\ =\ \frac{1}{2^{j+1}} F_{q+2}^r F_{q +4}^{r'}
\end{eqnarray}
with $q,r,r'$ as given in \eqref{eqn: k} and \eqref{eqn: m}. Arguing similarly leads to formulas for the other three cases, which we state below.

\begin{prop}\label{prop: exact}
Consider $i,j$ such that $i < j$.

For $i, j$ both odd:
\begin{equation}
\pp(i \mbox{ and }j \not \in A+A)\ =\ \frac{1}{2^{j+1}} F_{q+2}^r F_{q +4}^{r'}
\end{equation}
where
\begin{eqnarray}
q & \ = \ & 2\left\lceil \frac{i+1}{j-i} \right\rceil \nonumber\\
r &=& \frac{1}{2}\left((j-i)\left\lceil \frac{i+1}{j-i} \right\rceil - (i +1)\right) \nonumber\\
r' &=& \frac{1}{2}\left(j+1 -(j-i)\left\lceil \frac{i+1}{j-i} \right\rceil \right).
\end{eqnarray}

For $i$ even, $j$ odd:
\begin{equation}
\pp(i \mbox{ and }j \not \in A+A)\ =\ \frac{1}{2^{j+1}} F_o F_{q+2}^r F_{q +4}^{r'}
\end{equation}
where
\begin{eqnarray}
o &\ = \ & 2 \left \lceil \frac{i/2+1}{j-i} \right \rceil -1 \nonumber\\
q &=& 2\left\lceil \frac{i+1}{j-i} \right\rceil \nonumber\\
r &=& \frac{1}{2}\left((j-i-1)\left\lceil \frac{i+1}{j-i} \right\rceil - (i +1)+o\right) \nonumber\\
r' &=& \frac{1}{2}\left(j -(j-i-1)\left\lceil \frac{i+1}{j-i} \right\rceil - o \right).
\end{eqnarray}

For $i$ odd, $j$ even:
\begin{equation}
\pp(i \mbox{ and }j \not \in A+A)\ =\ \frac{1}{2^{j+1}} F_{o'+2} F_{q+2}^r F_{q +4}^{r'}
\end{equation}
where
\begin{eqnarray}
o' &\ = \ &2 \left \lceil \frac{j/2+1}{j-i} \right \rceil -2 \nonumber\\
q &=& 2\left\lceil \frac{i+1}{j-i} \right\rceil \nonumber\\
r &=& \frac{1}{2}\left((j-i-1)\left\lceil \frac{i+1}{j-i} \right\rceil - (i +1)+ o'\right) \nonumber\\
r' &=& \frac{1}{2}\left(j -(j-i-1)\left\lceil \frac{i+1}{j-i} \right\rceil - o' \right).
\end{eqnarray}

For $i, j$ both even:
\begin{equation}
\pp(i \mbox{ and }j \not \in A+A)\ =\ \frac{1}{2^{j+1}} F_o F_{o'}F_{q+2}^r F_{q +4}^{r'}
\end{equation}
where
\begin{eqnarray}
o &\ =\ & 2 \left \lceil \frac{i/2+1}{j-i} \right \rceil -1 \nonumber\\
o' &=& 2 \left \lceil \frac{j/2+1}{j-i} \right \rceil -2 \nonumber\\
q &=& 2\left\lceil \frac{i+1}{j-i} \right\rceil \nonumber\\
r &=& \frac{1}{2}\left((j-i-2)\left\lceil \frac{i+1}{j-i} \right\rceil - (i +1)+ o + o'\right) \nonumber\\
r' &=& \frac{1}{2}\left(j-1 -(j-i-2)\left\lceil \frac{i+1}{j-i} \right\rceil - o - o'\right).
\end{eqnarray}
\end{prop}

We conclude this section with some bounds on  $\pp(i \mbox{ and }j \not \in A+A)$. We have (Binet's formula)
\begin{equation}
F_n\ =\ \frac{1}{\sqrt{5}} (\phi^n - (-1/\phi)^n),
\end{equation}
where $\phi = (1+\sqrt{5})/2$ is the golden ratio. Therefore, for even $n$ we have
\begin{equation}
F_n\ \le\ \frac{1}{\sqrt{5}} \phi^n.
\end{equation}
Since $q+2$ and $q+4$ are always even, then for any $i,j$ both odd, we have
\begin{eqnarray}\label{eqn: ijdecay}
\pp(i \mbox{ and } j \not \in A+A) &=& \frac{1}{2^{j+1}} F_{q+2}^r F_{q+4}^{r'} \nonumber\\
&\le& \frac{1}{2^{j+1}} \left(\frac{\phi^{q+2}}{\sqrt{5}} \right)^r \left(\frac{\phi^{q+4}}{\sqrt{5}} \right)^{r'} \nonumber\\
&=& \frac{1}{2^{j+1}} \frac{ \phi^{(qr + (q+2)r') + (2r + 2r')} }{ 5^{(r + r')/2}}\nonumber\\
&=& \frac{1}{2^{j+1}} \frac{\phi^{j+1 + j-i} }{ 5^{(j-i)/4}}\nonumber\\
&=& \frac{\phi^{2j+1}}{2^{j+1} 5^{j/4}} \frac{5^{i/4}}{ \phi^i},
\end{eqnarray}
where the second to last equality comes from \eqref{eqn: system}. In fact, we can use Proposition \ref{prop: exact} to show
that \eqref{eqn: ijdecay} holds for all $i,j$ (slightly better constants hold for the other $i,j$).

If $i = k$ and $j = k +m$, where $m$ is fixed and $k$ goes to infinity, a lower bound similar to \eqref{eqn: ijdecay} also holds. First note that for even $n$
\begin{eqnarray}
F_n^r &=& \frac{1}{5^{r/2}}\left(\phi^n - \phi^{-n}\right)^r\nonumber\\
&=& \frac{1}{5^{r/2}}\phi^{nr}  \left( 1- \phi^{-2n}\right)^r\nonumber\\
&=& \frac{1}{5^{r/2}}\phi^{nr}\left( 1 - r(1-c)^{r-1} \phi^{-2n}\right)
\end{eqnarray}
for some $c$ such that $0 < c < 1/\phi^{2n}$ by Taylor expansion.
Therefore for odd $i,j$, we have
\begin{eqnarray}\label{eqn: ijdecay2}
\pp(i \mbox{ and } j \not \in A+A) &=& \frac{1}{2^{j+1}} F_{q+2}^r F_{q+4}^{r'} \nonumber\\
&\ge& \frac{1}{2^{j+1}} \frac{1}{5^{(q+2)/2}}\phi^{(q+2)r} (1- r\phi^{-2(q+2)}) \frac{1}{5^{(q+4)/2}}\phi^{(q+4)r'} (1- r'\phi^{-2(q+4)})\nonumber\\
&=& \frac{\phi^{2j+1}}{2^{j+1} 5^{j/4}} \frac{5^{i/4}}{ \phi^i}(1- r\phi^{-2(q+2)})(1- r'\phi^{-2(q+4)})\nonumber\\
&\ge& \frac{\phi^{2j+1}}{2^{j+1} 5^{j/4}} \frac{5^{i/4}}{ \phi^i}(1- (r+r')\phi^{-2(q+2)}) \nonumber\\
&\ge& \frac{\phi^{2j+1}}{2^{j+1} 5^{j/4}} \frac{5^{i/4}}{ \phi^i}(1- (j-i)\phi^{-i/ (j-i)}),
\end{eqnarray}
and similar formulas hold for the other parity cases. If $j/i \rightarrow 1$ not too slowly, then the remainder term on the right-hand-side of \eqref{eqn: ijdecay2} goes to $1$.
For example, if $i = k$ and $j = k +m$, then we have the following corollary by combining \eqref{eqn: ijdecay}  and \eqref{eqn: ijdecay2}.

\begin{cor}\label{cor: ij}
For any fixed $m$,
\begin{equation}
\pp(k \mbox{ and }k+m \not \in A+A)\ \sim\ \frac{\phi^{2(k+m)+1}}{2^{(k+m)+1} 5^{(k+m)/4}} \frac{5^{k/4}}{\phi^{k}}
= \frac{\phi^{k+1} }{2^{k+1}} \frac{\phi^{2m}}{2^m 5^{m/4}},
\end{equation}
as $k$ goes to infinity with $k, k +m$ are both odd. Similar asympotics hold for general $k,k+m$. If we ignore the constants related to $m$, we have
\begin{equation}
\pp(k \mbox{ {\rm and} }k+m \not \in A+A) \ =\ \Theta( (\phi/2)^k)
\end{equation}
as $k$ goes to infinity with any $k, k +m$.
\end{cor}

Note that since $\pp(i \mbox{ and } j  \not \in A+A)$ has exponential decay in $i,j$ as seen in \eqref{eqn: ijdecay}, then \eqref{eqn: variancesimple} converges as $n\rightarrow \infty$; that is
\begin{equation}\label{eqn: limvar}
\lim_{n\rightarrow \infty} \Var\left(M_{[0, n-1]}(A)\right) = 4\sum_{i< j}\pp(i \mbox{ and } j \not \in A+A) - 40
\end{equation}
exists and is finite. In particular, we know that the limit is an infinite sum of Fibonacci products. However, we could not find a closed form for this sum.
Nonetheless, because of the exponential decay in the terms in the sum, we can approximate the variance well.
In particular, note that the tail of the sum has exponential decay:
\begin{eqnarray}
\sum_{n \le i < j} P(i \mbox{ and } j \not \in A+A) &\le&
\sum_{n \le i < j} \frac{\phi}{2} \left(\frac{\phi^2}{2\cdot 5^{1/4}} \right)^j \left(\frac{5^{1/4}}{\phi}\right)^i \nonumber\\
 &\le&  \frac{\phi}{2}
\left(\sum_{n \le j } \left(\frac{\phi^2}{2\cdot 5^{1/4}} \right)^j  \right)
\left(\sum_{n \le i } \left(\frac{5^{1/4}}{\phi}\right)^i  \right) \nonumber\\
&\le&  \frac{\phi}{2} \left(\frac{1}{1- \phi^2/2\cdot 5^{1/4}}\right)
\left(\frac{1}{1- 5^{1/4}/\phi} \right)
\left(\frac{\phi^2}{2\cdot 5^{1/4}} \right)^n \left(\frac{5^{1/4}}{\phi}\right)^n\nonumber\\
&\le& 87\left(\frac{\phi}{2}\right)^n \le 87(0.81)^n. \label{eqn: errorterm}
\end{eqnarray}
Here we use that \eqref{eqn: ijdecay} holds for all $i,j$.
Using Mathematica to sum the first $300$ terms of \eqref{eqn: limvar}, whose exact form is given in Proposition \ref{prop: exact}, we get the following approximation for the variance:
\begin{equation}
\lim_{n\rightarrow \infty} \Var\left(M_{[0, n-1]}(A)\right)=  35.9658 + E,
\end{equation}
where $|E|< 10^{-4}$. The error term $E$ comes mostly from truncating the computation of the 300-term series given by Mathematica.
By \eqref{eqn: errorterm}, the error term from truncating the series at $n=300$ is less than $87(0.81)^{300} \sim 3\cdot 10^{-28}$, which is
much less than the Mathematica error term.
This proves Theorem \ref{thm: variance}.

%%%%%%%%%%%%%%%%%%%%%%%%%%%%%%%%%%%%%%%%%%%%%%%%%%%%%%%%%%%%%%%%%%%%%%%%%%%%%%%%%%%%%%%%%%%%%%
%%%%%%%%%%%%%%%%%%%%%%%%%%%%%%%%%%%%%%%%%%%%%%%%%%%%%%%%%%%%%%%%%%%%%%%%%%%%%%%%%%%%%%%%%%%%%%
\section{Exponential Bounds}\label{sec: bounds}

We now prove Theorem \ref{thm: bounds} and find exponential bounds for the distribution of $M_{[0, n-1]}(A)$.

\begin{proof}[Proof of Theorem \ref{thm: bounds}]
For the lower bound, we construct many $A$ such that $A+A$ is missing $k$ elements. First suppose that $k$ is even. Let the first $k/2$ non-negative integers not be in $A$. Then let the rest of the elements of $A$ be any subset $A'$ that fills in (so $A'+A'$ has no missing elements between its largest and smallest elements); that is $M_{n-k/2}(A') = 0$. By~\cite{MO}*{Proposition 8}, we can show that
\begin{equation}\label{eqn: lowerboundonfull}
\pp(M_{[0, n-1]}(A') = 0) > 1/2^{10}
\end{equation}
independent of $n$. If $L \subseteq [0, \ell-1]$ and $U\subseteq [n-u, n-1]$ are fixed, then their proposition says that
\begin{equation}
\pp([2\ell -1, 2n -2u-1] \subseteq A'+A' \  | \  A'\cap [0, \ell-1] = L, A'\cap [n-u, n-1] = U) > 1 - 6(2^{-|L|}+ 2^{-|U|}),
\end{equation}
independent of $n$.
Therefore,
\begin{eqnarray}
&&\pp([2\ell -1, 2n -2u-1] \subseteq A'+A' \mbox{ and } A'\cap [0, \ell-1] = L, A'\cap [n-u, n-1] = U) \nonumber\\
&& \ > \ (1 - 6(2^{-|L|}+ 2^{-|U|}))2^{-\ell}2^{-u}.
\end{eqnarray}
By letting $L = [0, \ell-1], U = [n-u, n-1]$ so the ends fill in, we get that
\begin{equation}
\pp(A'+A' = [0, 2n-2]) \  > \ (1 - 6(2^{-\ell}+ 2^{-u}))2^{-\ell}2^{-u}.
\end{equation}
Letting $\ell = u = 4$\  so that the first term in the product is positive, we get that
\begin{equation}
\pp(A'+A' = [0, 2n-2]) \  > \ (1 - 6(2^{-4}+ 2^{-4}))2^{-4}2^{-4}\ =\ 1/2^{10},
\end{equation}
independent of $n$, which gives us \eqref{eqn: lowerboundonfull}.

As $A = k/2 + A'$, we have $A+A = k + A'+A' = [k, 2n-2]$ and so $M_{[0, n-1]}(A) = k$ as seen by Figure \ref{fig: lower}.
\begin{figure}[h]
\includegraphics[scale=0.6]{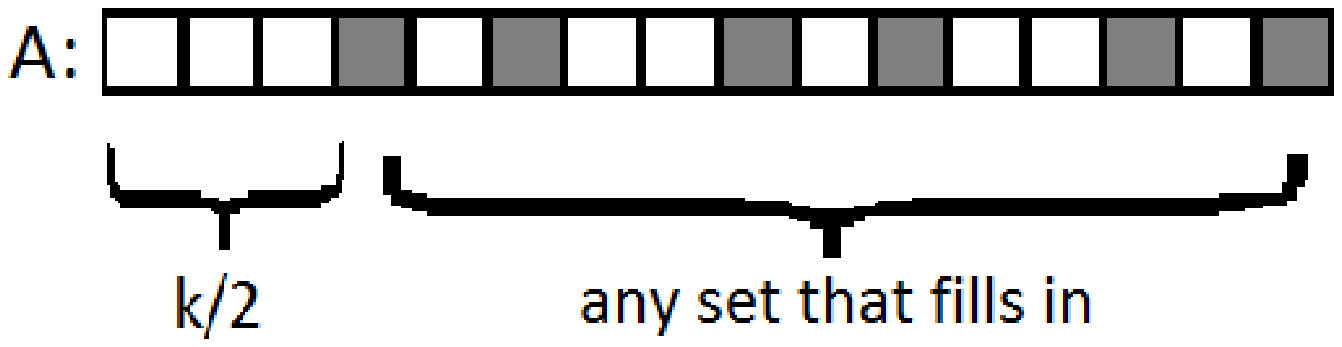}
\includegraphics[scale=0.6]{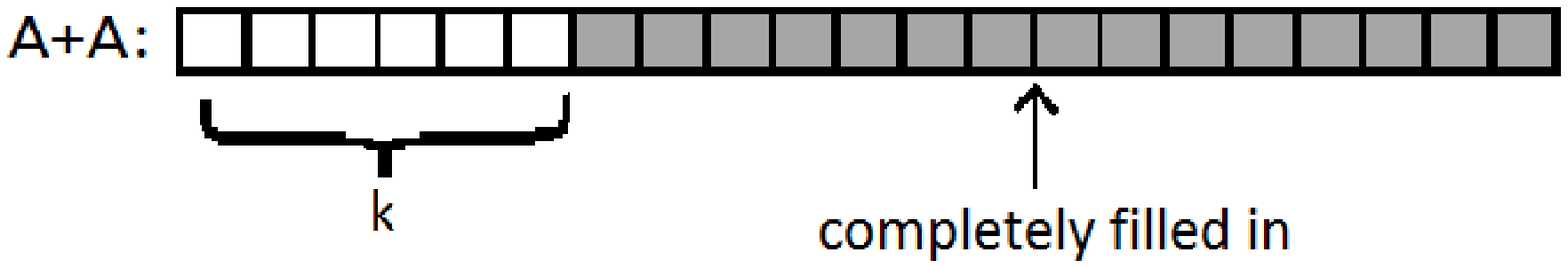}.
\caption{$A$ and $A+A$ for lower bound.}
\label{fig: lower}
\end{figure}

Therefore we have
\begin{eqnarray}
\pp(M_{[0, n-1]}(A) = k) &\ge& \pp(A = k/2 + A' \mbox{ and } M_{n-k/2}(A') = 0) \nonumber\\
&=& \left(\frac{1}{2}\right)^{k/2} \pp(M_{n-k/2}(A') = 0) \nonumber\\
&\gg& \left(\frac{1}{2}\right)^{k/2} \ \ge \ (0.70)^k,\label{eqn: lowerboundeven} \label{eqn: lowerboundodd}
\end{eqnarray}
where the implied constants are independent of $n$ by \eqref{eqn: lowerboundonfull}.
This proves the lower bound in Theorem \ref{thm: bounds} when $k$ is even.

If $k$ is odd, then we can let $L = [0, \ell-1] \setminus \{ 2,3 \}$ and $U = [n-u, n-1]$ so that only the element $3$ is missing from $A' + A'$. Then we get a bound for $\pp(M_{[0, n-1]}(A')=1)$ as in \eqref{eqn: lowerboundonfull}. Letting $A = (k-1)/2 + A'$, we get the desired lower bound in Theorem \ref{thm: bounds} for when $k$ is odd.
\\

For the upper bound, we can use bounds like
\begin{equation}\label{eqn: 3/4}
\pp(k\not \in A+A) \le \left(\frac{3}{4}\right)^{k/2}
\end{equation}
from \cite{MO}. Again, first suppose that $k$ is even. Note that if $A+A$ is missing $k$ elements, then one of these missing elements must be at least $k/2$ elements away from the ends of $[0,2n-2]$. That is, we have the following situation (see Figure \ref{fig: upperbound}).
\begin{figure}[h]
\includegraphics[scale=0.85]{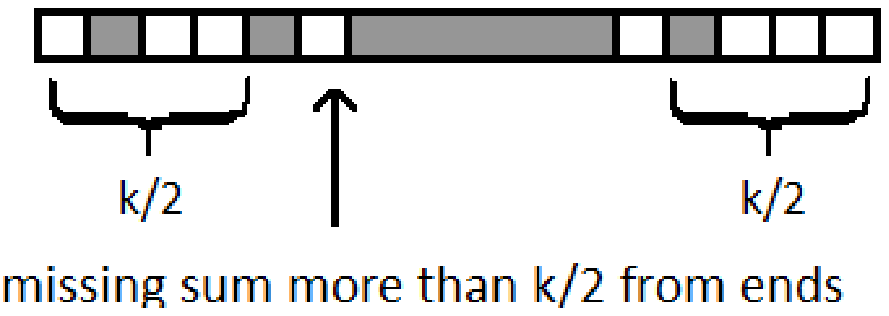}
\caption{Upper bound for $\pp(M_{[0, n-1]}(A)=k)$ }
\label{fig: upperbound}
\end{figure}

Therefore
\begin{eqnarray}
\pp(M_{[0, n-1]}(A) = k) &\le& \pp(A+A \mbox{ missing element at least } k/2 \mbox{ away from edges})\nonumber\\
&=\ & \pp(j \not \in A+A, j \in [k/2, 2n-k/2] )\nonumber\\
& \le\ & 2\sum_{j \ge k/2} \left(\frac{3}{4}\right)^{j/2} \nonumber\\
& \ll\ & \left(\frac{3}{4} \right)^{k/4} \approx (0.93)^k.
\end{eqnarray}

Note that this bound does not use the fact that there may be missing elements on both ends
at the same time. By focusing on one particular side, we can get a stronger result. For example, we have the
following inequality for the probability of missing $k$ elements in $[0,n/2]$:
\begin{eqnarray}
\pp( |[0,n/2] \setminus (A+A)| = k) & \le &  \pp(j \not \in A+A, j \in [k, n/2])\nonumber\\
&\ \le\ & 2\sum_{j \ge k} \left(\frac{3}{4}\right)^{j/2} \nonumber\\
&\ \ll\ & \left(\frac{3}{4} \right)^{k/2} \approx (0.87)^k \label{eqn: inequality 0.87}
\end{eqnarray}
and similarly for $\pp(|[3n/2,2n] \setminus (A+A)| = k).$
Furthermore, \eqref{ineq: appendix} from Section \ref{sec:bimodality}
connects the probability of missing $k$ elements to the probability of missing elements in
$[0,n/2]$ and $[3n/2, 2n]$:
\begin{equation}\label{eqn: eqappendix}
 \pp(M_{[0, n-1]}(A) = k) = \sum_{i+j = k} \pp( |[0,n/2] \setminus (A+A)| = i)\pp(|[3n/2,2n] \setminus (A+A)| = j )  + O\left(\left(\frac{3}{4}\right)^{n/4}\right).
\end{equation}
Combining \eqref{eqn: inequality 0.87} and \eqref{eqn: eqappendix}, we get
\begin{eqnarray}
&& \pp(M_{[0, n-1]}(A) = k) \nonumber\\
&&= \sum_{i+j = k} \pp( |[0,n/2] \setminus (A+A)| = i)\pp(|[3n/2,2n] \setminus (A+A)| = j )  + O\left(\left(\frac{3}{4}\right)^{n/4}\right)\nonumber\\
&&\ll\ \sum_{i+j = k} \left(\frac{3}{4} \right)^{i/2} \left(\frac{3}{4} \right)^{j/2}  + \left(\frac{3}{4}\right)^{n/4}\nonumber\\
&&\ll\ k \left(\frac{3}{4} \right)^{k/2}  + \left(\frac{3}{4}\right)^{n/4}.
\end{eqnarray}
Therefore if $k/2 < n/4$, then we get the desired bound
\begin{equation}\label{eqn: badupperbound}
\pp(M_{[0, n-1]}(A)= k)\ \ll\ k\left(\frac{3}{4}\right)^{k/2}\ \approx\ (0.87)^k.
\end{equation}

Note that the bound in \eqref{eqn: badupperbound} for the distribution is exactly the same as the bound in \eqref{eqn: inequality 0.87} for missing elements on a single side. Since all our bounds are exponential and \eqref{eqn: eqappendix} multiplies $\pp( |[0,n/2] \setminus (A+A)| = i)$ with
$\pp(|[3n/2,2n] \setminus (A+A)| = j )$, we can always use this approach to transform bounds on the probability of missing elements in $(A+A)\cap [0,n/2]$ to equally good bounds on number of missing elements in all of $A+A$. So it is sufficient to just develop bounds on missing elements on one side of $A+A$. In particular, we can use this approach to transform the bounds in Corollary \ref{cor: ij} to improve the bounds in \eqref{eqn: badupperbound}. By Corollary \ref{cor: ij}, we have
\begin{eqnarray}
\pp(|[0,n/2] \setminus (A+A)| = k) & \ \le \ &  \pp(A+A \mbox{ misses 2 elements greater than } k-3)\nonumber\\
& \ = \ &  \pp(i, j \not \in A+A, i, j \in [k-3, n/2])\nonumber\\
& = & \sum_{k-3 < i < j} \pp(i \mbox{ and } j \not \in A+A)\nonumber\\
& \ll &\sum_{k-3 < i < j}  \frac{\phi^{2j+1}}{2^{j+1} 5^{j/4}} \frac{5^{i/4}}{ \phi^i}\nonumber\\
& \ll &\frac{\phi^{2k+1}}{2^{k+1} 5^{k/4}} \frac{5^{k/4}}{ \phi^k} = \left(\frac{\phi}{2}\right)^k \approx (0.81)^k.
\end{eqnarray}
Then using the previous approach, we get a similar bound on the total number of missing sums:
\begin{equation}\label{eqn: improvedbound}
\pp(M_{[0, n-1]}(A) = k )\ \ll\ \left(\frac{\phi}{2}\right)^k \approx (0.81)^k.
\end{equation}

Note that as in \eqref{eqn: eqappendix}, we always have an extra $(3/4)^{n/4}$ term.
To make this term negligible, we need to have $(3/4)^{n/4} < (0.81)^k$, which means
$n > k \cdot 4\log (0.81)/ \log (3/4) \sim 2.92 k$ or that
$k < 0.34n$. This condition is sufficient in this case where we have the bound $(\phi/2)^k$.
However in general, we know that we have a lower bound of $(1/2)^{k/2}$ for the distribution. Therefore, to make the  $(3/4)^{n/4}$ term always negligible, we can have $(3/4)^{n/4} < (1/2)^{k/2}$, which means \\ $n > k \cdot 2\log(1/2)/ \log(3/4) \sim 5 k$, as in the statement of Theorem \ref{thm: bounds}. Note that then the implied constants are independent of $n$. Combining \eqref{eqn: lowerboundodd} and \eqref{eqn: improvedbound}, we get Theorem \ref{thm: bounds}.
\end{proof}

%%%%%%%%%%%%%%%%%%%%%%%%%%%%%%%%%%%%%%%%%%%%%%%%%%%%%%%%%%%%%%%%%%%%%%%%%%%%%%%%%%%%%%%%%%%%%%
%%%%%%%%%%%%%%%%%%%%%%%%%%%%%%%%%%%%%%%%%%%%%%%%%%%%%%%%%%%%%%%%%%%%%%%%%%%%%%%%%%%%%%%%%%%%%%
\section{Approximating $\pp (k+a_1, k+a_2, \dots, \mbox{ and } k+a_m \not \in A+A)$}\label{sec: configuration}

In this section, we prove Theorem \ref{thm: configuration}
which says that for any fixed $a_1, \dots, a_m$, there exists $\lambda_{a_1, \dots, a_m}$ such that
\begin{equation}
\pp(k + a_1, k + a_2, \dots, \mbox{ and } k + a_m \not \in A+A)\ =\ \Theta(\lambda_{a_1, \dots, a_m}^k),
\end{equation}
where the implied constants depend on $a_1, \dots, a_m$ but not $k$. Therefore, the probability is approximately exponential.

To prove this theorem, we use a version of Fekete's Lemma, which says that sub-additive sequences are approximately linear. From \cite{St} we have the following version in which the sequence is both sub-additive and super-additive.

\begin{lem}\label{lem: fekete}
If $b_n$ is a sequence such that
\begin{equation}
b_n + b_m -1\ \le\ b_{n+m}\ \le\ b_n + b_m +1
\end{equation}
for all $n,m$, then $\lambda = \inf b_n/n$ exists and for all $n$,
\begin{equation}
\left|\frac{b_n}{n} - \lambda \right|\  <\ \frac{1}{n}.
\end{equation}
\end{lem}

\begin{rek}
The proof of this Lemma can be easily modified to get that if
\begin{equation}
b_n + b_m -c\ \le\ b_{n+m}\ \le\ b_n + b_m +c
\end{equation}
for some constant $c>0$, then
\begin{equation}
\left|\frac{b_n}{n} - \lambda \right|\  <\ \frac{c}{n}.
\end{equation}
\end{rek}

Suppose that $a_n$ is approximately multiplicative rather than approximately additive so that for some constant $c >1$
\begin{eqnarray}\label{eqn: multiplicative}
c^{-1} \cdot a_m a_n\ \le\ a_{m+n}\ \le\ c\cdot a_m a_n
\end{eqnarray}
for all $m,n$. As $b_n = \log a_n$ satisfies the properties of Lemma \ref{lem: fekete}, for $\lambda = \inf \frac{\log a_n}{n}$ we have
\begin{equation}
\left|\frac{\log a_n}{n} - \lambda \right|\  <\ \frac{\log c}{n}
\end{equation}
for all $n$. That is,
\begin{equation}
c^{-1}\lambda^n\ \le\   a_n\ \le\ c \lambda^n
\end{equation}
for all $n$, implying
\begin{eqnarray}
a_n\ =\ \Theta(\lambda^n).
\end{eqnarray}
Therefore we just need to relate $\pp(k + a_1, k + a_2, \dots, \mbox{ and } k + a_m \not \in A+A)$ as a function of $k$ to some approximately multiplicative function satisfying \eqref{eqn: multiplicative}.

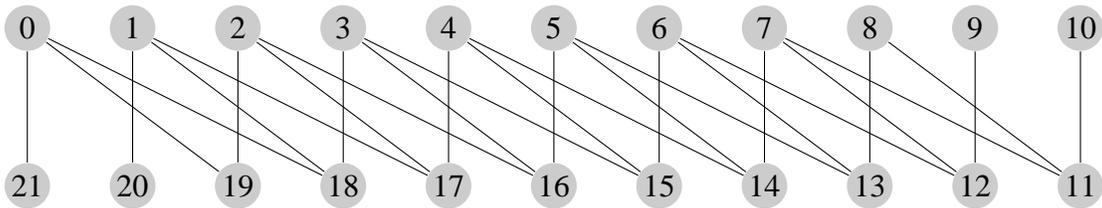
\begin{figure}[h]
\begin{tikzpicture}
   [scale=.7,auto=left,every node/.style={circle,fill=black!20,minimum size=2pt,inner sep=1pt}]
  \node (n0) at (0,1)  {\ 0 \ };
  \node (n1) at (2,1)  {\ 1 \ };
  \node (n2) at (4,1)  {\ 2 \ };
  \node (n3) at (6,1)  {\ 3 \ };
	\node (n4) at (8,1)  {\ 4 \ };
  \node (n5) at (10,1)  {\ 5 \ };
  \node (n6) at (12,1)  {\ 6 \ };
  \node (n7) at (14,1)  {\ 7 \ };
	\node (n8) at (16,1)  {\ 8 \ };
	\node (n9) at (18,1)  {\ 9 \ };
	\node (n10) at(20,1)  { 10 };
	
	\node (n21) at (0,-2)   {21};
  \node (n20) at (2,-2)   {20};
  \node (n19) at (4,-2)   {19};
  \node (n18) at (6,-2)   {18};
	\node (n17) at (8,-2)   {17};
  \node (n16) at (10,-2)  {16};
  \node (n15) at (12,-2)  {15};
  \node (n14) at (14,-2)  {14};
	\node (n13) at (16,-2)  {13};
	\node (n12) at (18,-2)  {12};
	\node (n11) at (20,-2)  {11};
	
  \foreach \from/\to in {
  											 n0/n21, n0/n19,n0/n18,
                         n1/n20, n1/n18,n1/n17,
                         n2/n19, n2/n17,n2/n16,
                         n3/n18, n3/n16,n3/n15,
                         n4/n17, n4/n15,n4/n14,
                         n5/n16, n5/n14,n5/n13,
                         n6/n15, n6/n13,n6/n12,
                         n7/n14, n7/n12,n7/n11,
                         n8/n13, n8/n11,
                         n9/n12,
                         n10/n11}
    \draw (\from) -- (\to);

\end{tikzpicture}
\caption{Modified condition graph for  $\pp(18,19, 21 \not \in A+A)$.}
\label{fig: modified}
\end{figure}

For example, consider $\pp(18,19, \mbox{ and } 21 \not \in A+A)$, whose condition graph is in Figure \ref{fig: condition161719}. Note that this graph has a loop from vertex 9 to itself since $9+9 =18$. We can symmetrize this graph by removing this loop and also removing the edge between vertices $8$ and $10$ and the edge between vertices $9$ and $10$, resulting in the modified condition graph in Figure \ref{fig: modified}.
\begin{figure}[h]
\begin{tikzpicture}
   [scale=.7,auto=left,every node/.style={circle,fill=black!20,minimum size=2pt,inner sep=1pt}]
  \node (n0) at (0,1)  {\ 0 \ };
  \node (n1) at (2,1)  {\ 1 \ };
  \node (n2) at (4,1)  {\ 2 \ };
  \node (n3) at (6,1)  {\ 3 \ };
	\node (n4) at (8,1)  {\ 4 \ };
  \node (n5) at (10,1)  {\ 5 \ };
  \node (n6) at (12,1)  {\ 6 \ };
  \node (n7) at (14,1)  {\ 7 \ };
	\node (n8) at (16,1)  {\ 8 \ };
	\node (n9) at (18,1)  {\ 9 \ };
	\node (n10) at(20,1)  { 10 };
	
	\node (n21) at (0,-2)   {21};
  \node (n20) at (2,-2)   {20};
  \node (n19) at (4,-2)   {19};
  \node (n18) at (6,-2)   {18};
	\node (n17) at (8,-2)   {17};
  \node (n16) at (10,-2)  {16};
  \node (n15) at (12,-2)  {15};
  \node (n14) at (14,-2)  {14};
	\node (n13) at (16,-2)  {13};
	\node (n12) at (18,-2)  {12};
	\node (n11) at (20,-2)  {11};
	
  \foreach \from/\to in {
  											 n0/n21, n0/n19,n0/n18,
                         n1/n20, n1/n18,n1/n17,
                         n2/n19, n2/n17,n2/n16,
                         n3/n18, n3/n16,n3/n15,
                         n4/n17, n4/n15,n4/n14,
                         n5/n16, n5/n14,n5/n13,
                         n6/n15, n6/n13,n6/n12,
                         n7/n14, n7/n12,n7/n11,
                         n8/n13, n8/n11,
                         n9/n12, n9/n10,
                         n10/n11}
    \draw (\from) -- (\to);
    \draw (n9) edge[loop] (n9);
     \draw (n8) edge[bend left] (n10);
\end{tikzpicture}
\caption{Condition graph for  $\pp(18,19, 21 \not \in A+A)$.}
\label{fig: condition161719}
\end{figure}

\emph{Denote the probability of getting a vertex cover for graphs like the one in Figure \ref{fig: modified} of length $n$ by $f(n)$}; so the probability of getting a vertex cover in Figure \ref{fig: modified} is $f(11)$.

Note that $f(11)$ is an upper bound for the probability in the original condition graph in Figure \ref{fig: condition161719} since we have removed some edges.  On the other hand, we have the following lower bound:
\begin{eqnarray}\label{eqn: conditionalprob}
&& \pp(18, 19, \mbox{ and } 21 \not \in A+A) \nonumber\\
&& \ \ge \  \pp(18, 19,  21 \not \in A+A \mbox{ and } 9, 10, 11, 12\not \in A) \nonumber\\
&& \ = \ \pp(18, 19,  21 \not \in A+A \mid 9, 10,11, 12 \not \in A)  \pp(9, 10, 11, 12\not \in A).
\end{eqnarray}
Note that the condition graph for $\pp(18, 19, 21 \not \in A+A \mid 9,10, 11,12 \not \in A)$ is the original condition graph in Figure \ref{fig: condition161719} with all edges incident on vertices $9,10,11$ or $12$ removed, as depicted in Figure \ref{fig: modified2}.

\begin{figure}[h]
\begin{tikzpicture}
   [scale=.7,auto=left,every node/.style={circle,fill=black!20,minimum size=2pt,inner sep=1pt}]
  \node (n0) at (0,1)  {\ 0 \ };
  \node (n1) at (2,1)  {\ 1 \ };
  \node (n2) at (4,1)  {\ 2 \ };
  \node (n3) at (6,1)  {\ 3 \ };
	\node (n4) at (8,1)  {\ 4 \ };
  \node (n5) at (10,1)  {\ 5 \ };
  \node (n6) at (12,1)  {\ 6 \ };
  \node (n7) at (14,1)  {\ 7 \ };
	\node (n8) at (16,1)  {\ 8 \ };
	
	\node (n21) at (0,-2)   {21};
  \node (n20) at (2,-2)   {20};
  \node (n19) at (4,-2)   {19};
  \node (n18) at (6,-2)   {18};
	\node (n17) at (8,-2)   {17};
  \node (n16) at (10,-2)  {16};
  \node (n15) at (12,-2)  {15};
  \node (n14) at (14,-2)  {14};
	\node (n13) at (16,-2)  {13};
	
  \foreach \from/\to in {
  											 n0/n21, n0/n19,n0/n18,
                         n1/n20, n1/n18,n1/n17,
                         n2/n19, n2/n17,n2/n16,
                         n3/n18, n3/n16,n3/n15,
                         n4/n17, n4/n15,n4/n14,
                         n5/n16, n5/n14,n5/n13,
                         n6/n15, n6/n13,
                         n7/n14,
                         n8/n13}
                         \draw (\from) -- (\to);

\end{tikzpicture}
\caption{Condition graph for $\pp(18, 19, \mbox{ and } 21 \not \in A+A \mid 9, 10, 11,12 \not \in A)$.}
\label{fig: modified2}
\end{figure}
Note that in Figure \ref{fig: modified2} we have removed vertices $9,10,11$ and $12$ completely since there are no longer any conditions on them in
$\pp(18, 19, \mbox{ and } 21 \not \in A+A \mid 9,10, 11,12 \not \in A)$. Finally, note that the probability of getting a vertex cover in the graph in Figure
\ref{fig: modified2} is just $f(9)$.
Therefore, by \eqref{eqn: conditionalprob}, we have
\begin{eqnarray}
 (1/2)^4 f(9)\ \le\ \pp(18,19,\mbox{ and } 21 \not \in A+A)\ \le\ f(11),
\end{eqnarray}
where we use that $\pp(9,10,11,12 \not \in A) = (1/2)^4$.

Since the condition graph for $\pp(k,k+1,\mbox{ and } k+3 \not \in A+A)$ is just a longer version of the condition graph for $\pp(18,19,\mbox{ and } 21 \not \in A+A)$, we can apply the same method as before to get that
\begin{eqnarray}\label{eqn: fp}
 (1/2)^4 f(k/2)\ \le\ \pp(k,k+1,\mbox{ and } k+3 \not \in A+A)\ \le\ f((k+4)/2)
\end{eqnarray}
for even $k$, with a similar formula holding for odd $k$. Therefore we are reduced to studying $f(n)$, which is easier to investigate since the condition graph is more symmetric. We will show that $f(n)$ satisfies \eqref{eqn: multiplicative}, implying it is approximately exponential.

\begin{figure}[h]
\begin{tikzpicture}
   [scale=.7,auto=left,every node/.style={circle,fill=black!20,minimum size=2pt,inner sep=1pt}]
  \node (n0) at (0,1)  {\ 0 \ };
  \node (n1) at (2,1)  {\ 1 \ };
  \node (n2) at (4,1)  {\ 2 \ };
  \node (n3) at (6,1)  {\ 3 \ };
	\node (n4) at (8,1)  {\ 4 \ };
  \node (n5) at (10,1)  {\ 5 \ };
  \node (n6) at (12,1)  {\ 6 \ };
  \node (n7) at (14,1)  {\ 7 \ };
	\node (n8) at (16,1)  {\ 8 \ };
	\node (n9) at (18,1)  {\ 9 \ };
	\node (n10) at(20,1)  { 10 };
	
	\node (n21) at (0,-2)   {21};
  \node (n20) at (2,-2)   {20};
  \node (n19) at (4,-2)   {19};
  \node (n18) at (6,-2)   {18};
	\node (n17) at (8,-2)   {17};
  \node (n16) at (10,-2)  {16};
  \node (n15) at (12,-2)  {15};
  \node (n14) at (14,-2)  {14};
	\node (n13) at (16,-2)  {13};
	\node (n12) at (18,-2)  {12};
	\node (n11) at (20,-2)  {11};
	
  \foreach \from/\to in {
  											 n0/n21, n0/n19,n0/n18,
                         n1/n20, n1/n18,
                         n2/n19,
                         n3/n18,
                         n4/n17, n4/n15,n4/n14,
                         n5/n16, n5/n14,n5/n13,
                         n6/n15, n6/n13,n6/n12,
                         n7/n14, n7/n12,n7/n11,
                         n8/n13, n8/n11,
                         n9/n12,
                         n10/n11}
    \draw (\from) -- (\to);
\end{tikzpicture}
\caption{Upper Bound for $f(11)$.}
\label{fig: upperf}
\end{figure}
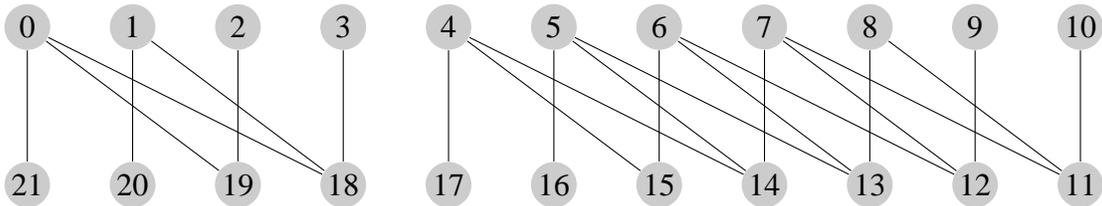

For example, to see that $f(11) \le f(4)f(7)$, we can separate the graph in Figure \ref{fig: modified} at the $4$\textsuperscript{th} vertex and remove edges that
cross this gap, resulting in the graph in Figure \ref{fig: upperf}.

Since the components are independent smaller copies of the original, the probability of getting a vertex cover for the graph in Figure \ref{fig: upperf} is $f(4) f(7)$. We can do this for any integer less than $11$, defining $f(n)$ for small integers by truncating at the $n$\textsuperscript{th} vertex. Since we have removed some edges to get the graph in Figure \ref{fig: upperf}, we have
\begin{eqnarray}
f(11) \ \le\ f(4)f(7)
\end{eqnarray}
as desired.

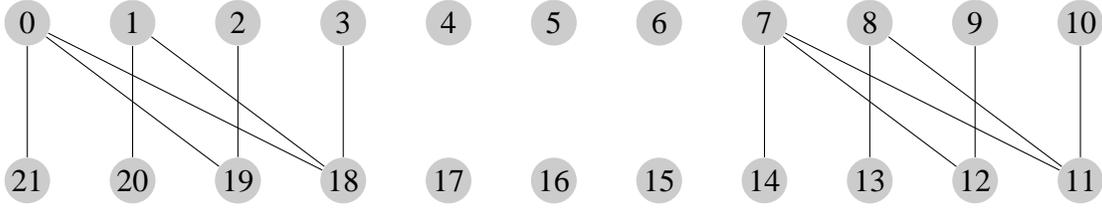
\begin{figure}[h]
\begin{tikzpicture}
     [scale=.7,auto=left,every node/.style={circle,fill=black!20,minimum size=2pt,inner sep=1pt}]
  \node (n0) at (0,1)  {\ 0 \ };
  \node (n1) at (2,1)  {\ 1 \ };
  \node (n2) at (4,1)  {\ 2 \ };
  \node (n3) at (6,1)  {\ 3 \ };
	\node (n4) at (8,1)  {\ 4 \ };
  \node (n5) at (10,1)  {\ 5 \ };
  \node (n6) at (12,1)  {\ 6 \ };  \node (n7) at (14,1)  {\ 7 \ };
	\node (n8) at (16,1)  {\ 8 \ };
	\node (n9) at (18,1)  {\ 9 \ };
	\node (n10) at(20,1)  { 10 };
	
	\node (n21) at (0,-2)   {21};
  \node (n20) at (2,-2)   {20};
  \node (n19) at (4,-2)   {19};
  \node (n18) at (6,-2)   {18};
	\node (n17) at (8,-2)   {17};
  \node (n16) at (10,-2)  {16};
  \node (n15) at (12,-2)  {15};
  \node (n14) at (14,-2)  {14};
	\node (n13) at (16,-2)  {13};
	\node (n12) at (18,-2)  {12};
	\node (n11) at (20,-2)  {11};
	
  \foreach \from/\to in {
  											 n0/n21, n0/n19,n0/n18,
                         n1/n20, n1/n18,
                         n2/n19,
                         n3/n18,
                         n7/n14, n7/n12,n7/n11,
                         n8/n13, n8/n11,
                         n9/n12,
                         n10/n11}
    \draw (\from) -- (\to);

\end{tikzpicture}
\caption{Lower Bound for $f(11)$.}
\label{fig: lowerf}
\end{figure}

To get a lower bound for $f(11)$, we use that
\begin{equation}\label{eqn: flower}
f(11)\ \ge\ f(11 \mid 4,5,6,15,16,17 \mbox{ chosen} ) \pp(4,5,6,15,16,17 \mbox{ chosen}),
\end{equation}
where $f(11 \mid 4,5,6,15,16,17 \mbox{ chosen})$ denotes the probability of getting a vertex cover for the graph in Figure \ref{fig: upperf} given that the vertices $4,5,6,15,16,17$ are chosen. The graph for $f(11 \mid 4,5,6,15,16,17 \mbox{ chosen })$ is depicted in Figure \ref{fig: lowerf}.
The probability of getting vertex covers for the two independent components is $f(4)f(4)$. Therefore from \eqref{eqn: flower}, we get that
\begin{equation}
f(11)\ \ge\ (1/2)^6 f(4)f(4)\ \ge\ (1/2)^6 f(4)f(7),
\end{equation}
with the last inequality since $f(n)$ is decreasing.
Therefore, in general we have
\begin{equation}
(1/2)^6 f(m)f(n)\ \le\ f(m+n)\ \le\ f(m)f(n),
\end{equation}
and so $f(n)$ satisfies the conditions of \eqref{eqn: multiplicative}. By the modified version of Fekete's Lemma, we have
\begin{equation}
f(n)\ =\ \Theta(\lambda^n)
\end{equation}
for some $\lambda$. Therefore by \eqref{eqn: fp}, we have
\begin{eqnarray}
\pp(k,k+1,\mbox{ and } k+3 \not \in A+A)\ =\ \Theta(\lambda^{k/2}),
\end{eqnarray}
which proves Theorem \ref{thm: configuration} for the case $a_1 = 0, a_2 = 1, a_3 = 3$.

The general situation follows in exactly the same way: by first making the configuration graph of $\pp(k+a_1, \dots, \mbox{ and } k + a_m \not \in A+A)$ look more symmetric and then using the modified Fekete's Lemma.

%%%%%%%%%%%%%%%%%%%%%%%%%%%%%%%%%%%%%%%%%%%%%%%%%%%%%%%%%%%%%%%%%%%%%%%%%%%%%%%%%%%%%%%%%%%%%%
%%%%%%%%%%%%%%%%%%%%%%%%%%%%%%%%%%%%%%%%%%%%%%%%%%%%%%%%%%%%%%%%%%%%%%%%%%%%%%%%%%%%%%%%%%%%%%
\section{Consecutive Missing Sums}\label{sec: consecutive}

In this section, we prove Theorem \ref{thm: consecutive} and its generalization Theorem \ref{thm: lambdabounds}.
We begin by proving Theorem \ref{thm: consecutive}, which says that
\begin{equation}
\left(\frac{1}{2}\right)^{(k+m)/2} \ll \pp(k+1, \dots, \mbox{ and } k + m \not \in A+A)\ \ll\ \left(\frac{1}{2}\right)^{(k+m)/2} (1+\epsilon_m)^k.
\end{equation}
The lower bound comes from the construction in Figure \ref{fig: lower} by letting the first $\lfloor (k+m)/2 \rfloor$ elements of $A$ be missing, which forces the first $k+m$ elements of $A+A$ to be missing as well.
That is,
\begin{eqnarray}\label{eqn: lowerboundexplained}
&& \pp(0, 1, \dots, k+m-1, \mbox{ and }k +m \not \in A+A)\nonumber\\
&& \: =\ \pp(0, 1, \dots, \mbox{ and }\lfloor (k+m)/2 \rfloor  \not \in A)\nonumber\\
&& \: = \ (1/2)^{\lfloor (k+m)/2 \rfloor +1}.
\end{eqnarray}
Therefore, we only need to prove the upper bound.

\begin{figure}[h]
\begin{tikzpicture}
   [scale=.7,auto=left,every node/.style={circle,fill=black!20,minimum size=2pt,inner sep=1pt}]
  \node (n0) at (0,1)  {\ 0 \ };
  \node (n1) at (2,1)  {\ 1 \ };
  \node (n2) at (4,1)  {\ 2 \ };
  \node (n3) at (6,1)  {\ 3 \ };
	\node (n4) at (8,1)  {\ 4 \ };
  \node (n5) at (10,1)  {\ 5 \ };
  \node (n6) at (12,1)  {\ 6 \ };
  \node (n7) at (14,1)  {\ 7 \ };
	\node (n8) at (16,1)  {\ 8 \ };
	
	\node (n20) at (-4,-2)  {20};	
	\node (n19) at (-2,-2)  {19};
  \node (n18) at (0,-2)   {18};
  \node (n17) at (2,-2)   {17};
  \node (n16) at (4,-2)   {16};
	\node (n15) at (6,-2)   {15};
  \node (n14) at (8,-2)   {14};
  \node (n13) at (10,-2)  {13};
  \node (n12) at (12,-2)  {12};
	\node (n11) at (14,-2)  {11};
	\node (n10) at (16,-2)  {10};
	\node (n9) at  (18,-2)  {\ 9 \ };
  \foreach \from/\to in {
  											 n0/n20, n0/n19,n0/n18,n0/n17,n0/n16,
                         n1/n19, n1/n18,n1/n17,n1/n16,n1/n15,
                         n2/n18, n2/n17,n2/n16,n2/n15,n2/n14,
                         n3/n17, n3/n16,n3/n15,n3/n14,n3/n13,
                         n4/n16, n4/n15,n4/n14,n4/n13,n4/n12,
                         n5/n15, n5/n14,n5/n13,n5/n12,n5/n11,
                         n6/n14, n6/n13,n6/n12,n6/n11,n6/n10,
                         n7/n13, n7/n12,n7/n11,n7/n10,n7/n9,
                         n8/n12, n8/n11,n8/n10,n8/n9,n8/n8,
                         n9/n10,n9/n9}
    \draw (\from) -- (\to);

     \draw (n9) edge[bend left] (n11);

\end{tikzpicture}
\caption{Condition graph for  $\pp(16,17,18,19,20 \not \in A+A)$.}
\label{fig:consecelementsone}
\end{figure}
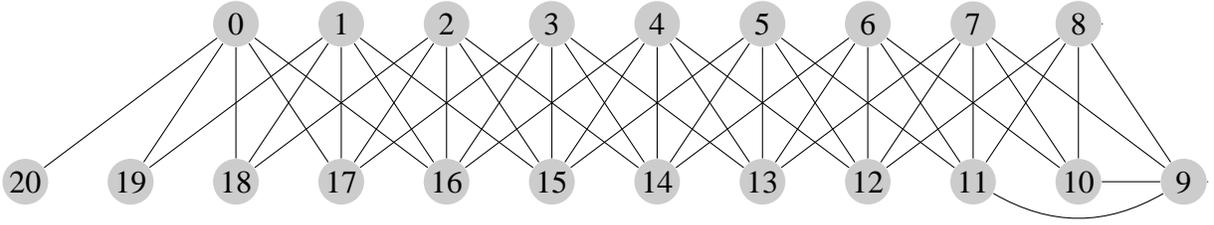

Before giving the proof, we consider an example with condition graphs which illustrates the idea. Consider $\pp(16,17,18,19,20 \not \in A+A)$. The condition graph here is given in Figure \ref{fig:consecelementsone}.

\begin{figure}[h]
\begin{tikzpicture}
   [scale=.7,auto=left,every node/.style={circle,fill=black!20,minimum size=2pt,inner sep=1pt}]
  \node (n0) at (0,1)  {\ 0 \ };
  \node (n1) at (2,1)  {\ 1 \ };
  \node (n2) at (4,1)  {\ 2 \ };
  \node (n3) at (6,1)  {\ 3 \ };
	\node (n4) at (8,1)  {\ 4 \ };
  \node (n5) at (10,1)  {\ 5 \ };
  \node (n6) at (12,1)  {\ 6 \ };
  \node (n7) at (14,1)  {\ 7 \ };
	\node (n8) at (16,1)  {\ 8 \ };
		
  \node (n20) at (-4,-2)  {20};		
	\node (n19) at (-2,-2)  {19};
  \node (n18) at (0,-2)   {18};
  \node (n17) at (2,-2)   {17};
  \node (n16) at (4,-2)   {16};
	\node (n15) at (6,-2)   {15};
  \node (n14) at (8,-2)   {14};
  \node (n13) at (10,-2)  {13};
  \node (n12) at (12,-2)  {12};
	\node (n11) at (14,-2)  {11};
	\node (n10) at (16,-2)  {10};
	\node (n9) at  (18,-2)  {\ 9 \ };
  \foreach \from/\to in {n0/n18,n0/n17,n0/n16,
                         n1/n18,n1/n17,n1/n16,
                         n2/n18,n2/n17,n2/n16,
                         n3/n15,n3/n14,n3/n13,
                         n4/n15,n4/n14,n4/n13,
                         n5/n15,n5/n14,n5/n13,
                         n6/n12,n6/n11,n6/n10,
                         n7/n12,n7/n11,n7/n10,
                         n8/n12,n8/n11,n8/n10}
    \draw (\from) -- (\to);
\end{tikzpicture}
\caption{Graph after removing some edges.}
\label{fig:figafterremovingsome}
\end{figure}
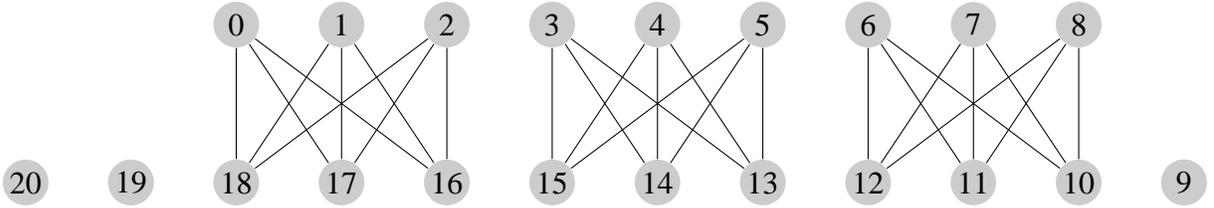

We need to find the probability of getting a vertex cover for this graph. If we remove some edges, the probability of getting a vertex cover for the resulting graph is an upper bound for the probability of getting a vertex cover for the original graph. We can remove some edges to get the graph of Figure \ref{fig:figafterremovingsome}.

The resulting graph has $3 \sim 20/6$ components that are all complete bipartite graphs with $6$ vertices. These are easier to handle since the only way to get a vertex cover for such graphs is to have all vertices on one side be chosen. So the probability of getting a vertex cover for one of these complete bipartite components is less than $(1/2)^3+ (1/2)^3 = 2/2^3$. Since the components are also independent, we have
\begin{equation}
\pp(16,17,18,19,20\not \in A+A)\ \le\ \left(\frac{2}{2^3}\right)^{3}\ \sim\ \left(\frac{1}{4}\right)^{20/6}.
\end{equation}
and in general we get that
\begin{equation}\label{eqn: condition01234}
\pp(k,k+1,k+2,k+3,k+4 \not \in A+A)\ \le\ \left(\frac{2}{2^3}\right)^{(k+4)/6}\  = \ \left(\frac{2^{1/3}}{2}\right)^{(k+4)/2}.
\end{equation}

We use this approach in the general proof. Notice that as $m \rightarrow \infty$, the size of the complete bipartite graphs          grows, and so we will be taking out relatively fewer and fewer constraints. Therefore, this approach gets us closer to the correct answer.

Now we give a formal proof of Theorem \ref{thm: consecutive} that does not rely on the condition graphs.
\begin{proof}
We first do the proof for $\pp(k, k +1, \dots, \mbox{ and } k + 2m-1 \not \in A+A)$ with $2m-1$ instead of $m$. Note that since the probability
depends only on $[0, k+2m-1] \cap A$, we can assume that $A \subseteq [0, k + 2m-1]$.
We will also assume that $m$ divides $k$ and that
\begin{equation}
k = qm
\end{equation}
with $q$ even.

We begin by writing $A$ as the following disjoint union:
\begin{equation}
A\ =\ A_0 \cup A_1 \cup \cdots \cup A_q \cup A_{q+1},
\end{equation}
where
\begin{equation}
A_j\ =\ A \cap [jm,(j+1)m-1].
\end{equation}
Then if $ [k, k+2m-1]\cap (A+A) = \emptyset$, then $[k, k+2m-1] \cap(A_j + A_{q-j}) = \emptyset$ for all $j$.
Note that
\begin{equation}
A_j + A_{q-j}\ \subseteq\ [k, k+2m -2].
\end{equation}
Therefore, $[k, k+2m-1] \cap(A_j + A_{q-j}) = \emptyset$ implies $A_j + A_{q-j} = \emptyset$. If $q$ is even, we have
\begin{eqnarray}\label{eqn: help1}
\pp(k, k+1,\dots, \mbox{ and } k + 2m -1 \not \in A+A)
&<& \pp([k, k+2m] \cap (A_j+A_{q-j}) = \emptyset \mbox{ for all } j\le q/2)\nonumber\\
&=& \pp(A_j+A_{q-j} = \emptyset \mbox{ for all } j\le q/2)\nonumber\\
&=& \pp(A_j = \emptyset \mbox{ or } A_{q-j} = \emptyset \mbox{ for all } j\le q/2).
\end{eqnarray}
For different $j$, the pairs of sets $A_j, A_{q-j}$ are disjoint. Therefore, we have independence:
\begin{equation}\label{eqn: help2}
\pp(A_j = \emptyset \mbox{ or } A_{q-j} = \emptyset \mbox{ for all } j\le q/2)
\ =\ \pp(A_{q/2} = \emptyset)\prod_{j = 0}^{q/2-1} \pp(A_j = \emptyset \mbox{ or } A_{q-j} = \emptyset).
\end{equation}
Finally, note that
\begin{equation}\label{eqn: help3}
\pp(A_j = \emptyset \mbox{ or } A_{q-j} = \emptyset)\ \le\ \pp(A_j = \emptyset) + \pp( A_{q-j} = \emptyset) = \frac{2}{2^m}.
\end{equation}
Combining \eqref{eqn: help1}, \eqref{eqn: help2}, and \eqref{eqn: help3}, we find
\begin{eqnarray}
\pp(k, k+1,\dots, \mbox{ and } k + 2m -1 \not \in A+A)
&\le& \frac{1}{2^m} \prod_{j=0}^{q/2-1} \frac{2}{2^m}\nonumber\\
&=&   2^{q/2}\left(\frac{1}{2^m}\right)^{q/2+1} \nonumber\\
&=&   2^{k/2m} \left(\frac{1}{2}\right)^{(k+2m)/2}.\label{eqn: }
\end{eqnarray}
This inequality is true for all $m,k$ such that $q = k/m$ is an even integer.

Changing $m$ to $m/2$, we get that
\begin{equation}\label{eqn: fin}
\pp(  k, k+1,\dots, \mbox{ and } k + m -1 \not \in A+A )\ \le\
2^{k/m} \left(\frac{1}{2}\right)^{(k+m)/2}
\end{equation}
for even $m$ and $q = k/m$ still an even integer. Note that \eqref{eqn: fin} is similar to the bound we get in \eqref{eqn: condition01234}   using the condition graph approach.

For odd $m$, we just need to use \eqref{eqn: fin}, noting that
\begin{equation}\label{eqn: fin2}
\pp( k, k+1,\dots, k + m -1, \mbox{ and } k +m \not \in A+A )\ \le\ \pp( k, k+1,\dots, \mbox{ and } k + m -1 \not \in A+A ).
\end{equation}
For odd $q$, we need to partition $A$ such that there is a block in the very middle of $A$. This ensures that this middle block is matched with itself (just like $A_{q/2}$ was matched with itself when $q$ was even). This gives us the extra $1/2^{m}$ that is needed in order to achieve the bound. For non-integer $q$, we need to repartition $A$ in a similar way. Therefore the bound in \eqref{eqn: fin} holds in general, up to a constant.

Finally, note that as $m\rightarrow \infty$, we have $2^{1/m} \rightarrow 1$. Writing $2^{1/m} = 1+ \epsilon_m$, we have
\begin{eqnarray}
\pp(k, \dots, \mbox{ and } k + m-1 \not \in A+A )\ <\ \left(\frac{1}{2}\right)^{(k+m)/2} (1+\epsilon_m)^{k},
\end{eqnarray}
where $\epsilon_m\rightarrow 0$ as $m\rightarrow \infty$. By raising $2^{1/m}= 1+ \epsilon_m$ to the $m$\textsuperscript{th} power, we see that
\begin{equation}
\epsilon_m < \frac{1}{m}.
\end{equation}
Therefore a weakened version of the inequality says that
\begin{equation}\label{eqn: weakened}
\left(\frac{1}{2}\right)^{(k+m)/2}\ \ll\ \pp(k+1, \dots, \mbox{ and } k +m \not \in A+A )\ \ll\ \left(\frac{1}{2}\right)^{(k+m)/2} (1+\epsilon_m)^{k},\end{equation}
where the implied constants are independent of $m$ and $k$.

This bound is interesting since it means that the trivial lower bound is almost the right answer for the exact bound. The trivial lower bound
makes us miss all of $[0, k+m]$ in $A+A$ as seen in \eqref{eqn: lowerboundexplained} but we only need $[k+1,k+m]$ to be missing.
In this sense, we see that essentially the only way to miss $m$ consecutive elements at $k+1$ for large $m$ is to miss all the previous elements as well.

Also, note that \eqref{eqn: weakened} implies that
\begin{equation}
\lambda_{0,1, \dots, m} \rightarrow \left(\frac{1}{2}\right)^{1/2}
\end{equation}
as $m\rightarrow \infty$ by definition of $\lambda_{0,1, \dots, m}$.
\end{proof}

Now we will prove Theorem \ref{thm: lambdabounds}, which says that
\begin{equation}\label{eqn: lambdabounds}
\lambda_{a_1, \dots, a_m}\ \le\ \pp(A,B \subseteq [0, \lfloor a_m/2 \rfloor ] \mid a_1, \dots, a_m \not \in A+B)^{1/(a_m+2)}.
\end{equation}
Note that Theorem \ref{thm: consecutive} is indeed a special case of this theorem since
we have the following upper bound
\begin{eqnarray}
\lambda_{0,1, \dots, m} &\le& \pp(A,B \subseteq [0, \lfloor m/2\rfloor ] \mid 0, \dots, m \not \in A+B)^{1/(m+2)}\nonumber\\
&=& \pp(A,B \subseteq [0, \lfloor m/2\rfloor ] \mid A = \emptyset \mbox{ or } B  = \emptyset)^{1/(m+2)}\nonumber\\
&\le& \left(2\left(\frac{1}{2}\right)^{\lfloor m/2\rfloor +1 }\right)^{1 /(m+2)}, % \nonumber\\
%&\rightarrow&  \left(\frac{1}{2}\right)^{1/2}.
\end{eqnarray} which converges to $\sqrt{1/2}$.

The proof of Theorem \ref{thm: lambdabounds} will be almost exactly the same as the proof of Theorem \ref{thm: consecutive}.

\begin{proof}
We will first show that for $0 \le a_1 < \cdots < a_m$,
\begin{eqnarray}
&&\pp(A \subseteq [0, k+a_m] \mid  k+a_1, \dots, k+a_m \not \in A+A) \nonumber\\
&&\:\:\:\:\:\:\le\ \pp(A,B \subseteq [0, a_m/2] \mid a_1, \dots, a_m \not \in A+B)^{k/(a_m+2)} \label{eqn: lambdainequality1}
\end{eqnarray}
for all $k,a_m$ such that $a_m$ is even and  $a_m+2$ divides $k$. Similar results hold in
the other cases of $k,m$. Furthermore, we first assume that $a_m = 2r-2$. Note that since the probability
depends only on $[0, k+2r-2] \cap A$, we can take $A \subseteq [0, k + 2r-2]$.
Again, we first assume that $r$ divides $k$ and that $k = qr$.
Then as before,
\begin{eqnarray}
&&\pp(k+a_1,\dots, \mbox{ and } k + a_m \not \in A+A)\nonumber\\
&&\: \: \: \: \: \: \le\ \pp(k+a_1,\dots, \mbox{ and } k + a_m \not \in A_j+A_{q-j} \mbox{ for all } j\le \lfloor q/2\rfloor)\nonumber\nonumber\\
&&\: \: \: \: \: \:  =\ \prod_{j=0}^{\lfloor q/2 \rfloor  }\pp(k+a_1,\dots, \mbox{ and } k + a_m \not \in A_j+A_{q-j} ).
\end{eqnarray}
The key fact is that if $j \ne q-j$, the sets $A_j, A_{q-j}$ are independent and
\begin{equation}
\pp(k+a_1,\dots, \mbox{ and } k + a_m \not \in A_j+A_{q-j} )\ =\ \pp(A, B \subseteq [0, r-1] \mid a_1, \dots, a_m \not \in A+ B)
\end{equation}
for all $j$. Therefore, if $q$ is odd
\begin{eqnarray}
&& \pp(k+a_1,\dots, \mbox{ and } k + a_m \not \in A+A) \nonumber\\
&& \: \: \: \: \: \:  \le\ \pp(A, B \subseteq [0, r-1] \mid a_1, \dots, a_m \not \in A+ B)^{\lfloor q/2\rfloor+1 }\nonumber\\
&&   \: \: \: \: \: \: =\ \pp(A, B \subseteq [0, r-1] \mid a_1, \dots, a_m \not \in A+ B)^{k/2r + 1/2}
\end{eqnarray}
and if $q$ is even,
\begin{eqnarray}
&& \pp(k+a_1,\dots, \mbox{ and } k + a_m \not \in A+A) \nonumber\\
&& \: \: \: \: \: \:  \le\  \pp(A \subseteq [0,r-1] \mid a_1, \dots, a_m \not \in A+A) \nonumber\\
&&\: \: \: \: \: \:\ \ \ \ \ \ \ \times\ \pp(A, B \subseteq [0, r-1] \mid a_1, \dots, a_m \not \in A+ B)^{k/2r}.
\end{eqnarray}
If we drop the terms that do not depend on $k$, we have for all even $a_m$ and all
$k$ divisible by $a_m+2$
\begin{eqnarray}
&&\pp(A \subseteq [0, k+a_m] \mid  k+a_1, \dots, k+a_m \not \in A+A) \nonumber\\
&&\:\:\:\:\:\:\le\ \pp(A,B \subseteq [0, a_m/2] \mid a_1, \dots, a_m \not \in A+B)^{k/(a_m+2)},
\end{eqnarray}
which is \eqref{eqn: lambdainequality1}.
Note that if $k$ is not divisible by $a_m+2$ or if $a_m$ is not even, we have
\begin{eqnarray}
&&\pp(A \subseteq [0, k+a_m] \mid  k+a_1, \dots, k+a_m \not \in A+A) \nonumber\\
&&\:\:\:\:\:\:\le\ \pp(A,B \subseteq [0,\lfloor a_m/2 \rfloor] \mid a_1, \dots, a_m \not \in A+B)^{\lfloor k/(a_m+2) \rfloor},
\end{eqnarray}
which proves that \eqref{eqn: lambdabounds}.

\end{proof}

%%%%%%%%%%%%%%%%%%%%%%%%%%%%%%%%%%%%%%%%%%%%%%%%%%%%%%%%%%%%%%%%%%%%%%%%%%%%%%%%%%%%%%%%%%%%%%%%%%%%%%%%%%%%%%%%
%%%%%%%%%%%%%%%%%%%%%%%%%%%%%%%%%%%%%%%%%%%%%%%%%%%%%%%%%%%%%%%%%%%%%%%%%%%%%%%%%%%%%%%%%%%%%%%%%%%%%%%%%%%%%%%%
%%%%%%%%%%%%%%%%%%%%%%%%%%%%%%%%%%%%%%%%%%%%%%%%%%%%%%%%%%%%%%%%%%%%%%%%%%%%%%%%%%%%%%%%%%%%%%%%%%%%%%%%%%%%%%%%

\section{Bounds on $m(k)$, $w(k)$, $y(k)$, and $z(k)$ for $k<32$}\label{sec:bimodality}

As mentioned in \S\ref{subsec: intro divot} and covered in more detail in \S\ref{subsec: connecting the distribs}, it suffices to bound $z(k)$. Our strategy is this: if $D+D$ (where $D$ is a uniformly chosen subset of \NN\ that contains 0) is missing exactly 7 elements, then it is very likely that those 7 missing sums are all smaller than 88 and typically even all smaller than 44. If we loop over all $2^{43}$ possibilities $\beta$ for $D\cap[0,44)$, for each possibility we can compute $(D+D)\cap[0,44)=(\beta+\beta)\cap[0,44)$ and a subset of $(D+D)\cap[44,48)\supseteq (\beta+\beta)\cap[44,88)$. From this (with some theory to handle the tail of the sumset) we bound the likelihood of missing exactly 7 sums, given $D\cap[0,44)$. By combining these estimates, we acquire bounds on $z(7)$.

Let $n\geq 2$ be a natural number parameter (the computations reported here use $n=44$, although $n=43$ is already enough to show $m(7)<m(6)<m(8)$), and set
   \be
      z(k \mid \beta)\ :=\ \Prob{|\NN \setminus (D+D)| = k \mid D \cap [0,n) = \beta}.
   \ee
We have
   \be
      z(k)\ =\ \sum_{0\in\beta\subseteq[0,n)} z(k\mid \beta) \Prob{D\cap[0,n) = \beta}\ =\ 2^{-(n-1)} \sum_{0\in\beta\subseteq[0,n)} z(k\mid \beta),
   \ee
so that it suffices to bound $z(k\mid\beta)$ above and below for all $0\leq k < 32$ (our arbitrary notion of ``small $k$'' is $0\le k < 32$) and all $0\in\beta\subseteq[0,n)$.

\newcommand{\Beginning}{{\mathcal B}}
\newcommand{\Definite}{{\mathcal D}}
\newcommand{\Likely}{{\mathcal L}}
\newcommand{\Tail}{{\mathcal T}}
\newcommand{\ExpectedMisses}{\eta}

Further, set
    \begin{align}
    \Beginning &\ :=\  D \cap [0,n) \nonumber\\
    \Definite &\ :=\  [0,n) \setminus (\beta+\beta) \nonumber\\
    \Likely &\ :=\  [n,2n) \setminus (\beta+\beta) \nonumber\\
    m &\ :=\  \min \Likely \nonumber\\
    \Tail &\ :=\  [2n,\infty) \nonumber\\
%    z(k\mid\beta) &\ :=\  \Prob{ | \NN \setminus (D+D) | = k \mid \Beginning=\beta } \nonumber\\
%    z(k) &\ :=\  \Prob{  | \NN \setminus (D+D) | = k } \nonumber\\
    \ExpectedMisses &\ :=\  \Expect{ | [n,\infty) \setminus (D+D) \mid \Beginning=\beta} \nonumber\\
    \mu &\ :=\  2^{-|\beta\cap [0,m-n]|}.
    \end{align}
If we condition on $\Beginning=\beta$, then the elements of $\Definite$ are $\Definite$efinitely missing from $D+D$, the elements of $\Likely$ are $\Likely$ikely but not certain to be missing, and the elements of $\Tail$, the $\Tail$ail of the natural numbers, are very likely to be missing. Note that $2n-1\in \Likely$, so $\Likely$ is nonempty and $m$ is well-defined.

\begin{lem}
   For all $k< |\Definite|$, we have $z(k \mid \beta) = 0$.
\end{lem}

\begin{proof}
   Conditioning on $\Beginning=\beta$, we have $\Definite\subseteq \NN \setminus(D+D)$. In fact, $\Definite=[0,n)\setminus(D+D)$.
\end{proof}

\begin{lem} We have
    \(\displaystyle
    \ExpectedMisses = 5 \cdot 2^{-|\Beginning|}+\sum_{\ell \in \Likely} 2^{-|\Beginning\cap[0,\ell-n]|}.
    \)
\end{lem}

\begin{proof}
By linearity of expectation
    \be
    \ExpectedMisses\ :=\ \Expect{ | [n,\infty) \setminus (D+D) | } = \Expect{ |[n,2n)\setminus(D+D)|} +  \Expect{ | \Tail \setminus(D+D)|}.
    \ee
Again using linearity of expectation, we have
    \be
    \Expect{ |[n,2n)\setminus(D+D)|}\ =\ \sum_{\ell \in \Likely} \Prob{\ell\not\in D+D}  %= \sum_{\ell \in \Likely} 2^{-|\Beginning\cap[0,\ell-n]|}
    \ee
Since $\ell\not\in D+D$ is the same as (for $n\leq \ell<2n$)
   \begin{align}
      \ell\not\in D+D
         \ =\  \bigwedge_{i=0}^{\ell/2} (i\not\in D \vee \ell-i\not\in D)
         \ =\ \bigwedge_{\substack{b\in \beta \\ b \leq \ell-n}} \ell-b \not\in D.
   \end{align}
Thus
   \be \Prob{\ell\not\in D+D}\ =\ 2^{-|\beta \cap [0,\ell-n]|}, \ee
and so
   \be \sum_{\ell \in \Likely} \Prob{\ell\not\in D+D}\ = \ \sum_{\ell \in \Likely} 2^{-|\beta\cap[0,\ell-n]|}.\ee

That
    \be \Expect{ | \Tail \setminus(D+D)|}\ =\  5 \cdot 2^{-|\beta|} \ee
is essentially in \cite{MO}, but we derive it here for the reader's convenience. By linearity of expectation,
   \begin{equation}\label{equ:T}
      \Expect{ | \Tail \setminus(D+D)|}\ =\ \sum_{t=2n}^\infty \Prob{t\not\in D+D}
   \end{equation}
and
   \begin{align}
      \Prob{t\not\in D+D}
            \ =\  \Prob{\left(\bigwedge_{{b\in \beta}} t-b\not\in D \right)
                  \bigwedge
               \left( \bigwedge_{i=n}^{t/2} i\not\in D \vee t-i \not\in D \right)}.
   \end{align}
Now this has two cases leading to
   \be
      \Prob{t\not\in D+D}\ =\ \begin{cases} 2^{-|\beta|} (3/4)^{(t-2n+1)/2} & \text{$t$ is odd,} \\ 2^{-|\beta|} (1/2) (3/4)^{(t-2n)/2} & \text{$t$ is even.} \end{cases}
   \ee
The infinite sum~\eqref{equ:T} now simplifies $5\cdot 2^{-|\beta|}$.
\end{proof}

\begin{lem}\label{lemma:z_D} We have
$\displaystyle \max\{0,1-\ExpectedMisses\} \leq z({|\Definite|} \mid \beta) \leq 1-\mu.$
\end{lem}

\begin{proof}
Trivially $z({|\Definite|} \mid \beta)\geq 0$. Since
\bea      \ExpectedMisses  & \ = \ & \Expect{|[n,\infty)\setminus (D+D)| \mid \Beginning=\beta}  \ =\ \sum_{i=0}^\infty z({|\Definite|+i} \mid \beta) \; \cdot \; i \nonumber\\
         &\ \geq\ & \sum_{i=1}^\infty z( {|\Definite|+i} \mid \beta) \ =\  1 - z({|\Definite|} \mid \beta),
\eea
we also have $z({|\Definite|} \mid \beta) \geq 1-\ExpectedMisses$.

Observe that the event $|\NN\setminus(D+D)| > |\Definite|$ contains the event $\{m\not\in D+D\}$, and so
   \be \Prob{|\NN\setminus(D+D)| = |\Definite|}\ \leq\ 1-\Prob{m\not\in D+D} = 1-\mu,\ee
concluding the proof of this lemma.
\end{proof}

\begin{lem}\label{lemma:z_D+1} We have
   $\displaystyle \max\{0,2\mu-\ExpectedMisses\} \leq z( {|\Definite|+1} \mid \beta) \leq \min\{1,\ExpectedMisses\}.$
\end{lem}

\begin{proof}
Trivially $0 \leq z({|\Definite|+1} \mid \beta) \leq 1$. We have
   \be
      \ExpectedMisses \ = \ \sum_{k=0}^\infty k\,\cdot\, z({|\Definite|+k} \mid \beta) \geq z({|\Definite|+1} \mid \beta),
   \ee
which leaves only the bound $2\mu-\ExpectedMisses \leq z({|\Definite|+1} \mid \beta)$ to prove.

The idea here is that if exactly $|\Definite|+1$ sums are missing, they are very likely to be the $|\Definite|$ elements of $\Definite$, and $m$. Formally,
   \begin{align*}
      \left\{ |\NN \setminus (D+D)|=|\Definite|+1\right\}
         & \supseteq \left\{ m \not\in D+D\right\} \cap \bigcap_{\substack{\ell\in \Likely \\ \ell>m}} \left\{ \ell\in D+D\right\} \cap \bigcap_{t\in \Tail} \left\{ t \in D+D\right\}\\
         & \supseteq \left\{ m \not\in D+D\right\} \setminus
               \left( \bigcup_{\substack{\ell\in \Likely \\ \ell>m}} \left\{ \ell\not\in D+D\right\} \cup \bigcup_{t\in \Tail} \left\{ t \not\in D+D\right\} \right)
   \end{align*}
and so
   \begin{align}
      z({|\Definite|+1} \mid \beta) &\ \geq\
         \Prob{m\not\in D+D} - \sum_{\substack{\ell\in \Likely \\ \ell>m}} \Prob{\ell \not\in D+D} - \sum_{t\in \Tail} \Prob{t\not\in D+D} \nonumber\\
         &\ = \ 2\Prob{m\not\in D+D} - \sum_{i\in \Likely \cup \Tail} \Prob{i \not\in D+D} \ = \ 2\mu-\ExpectedMisses.
   \end{align}
\end{proof}

\begin{lem}
For $k\geq 2$, $\displaystyle 0\leq z({|\Definite|+k}\mid \beta) \leq \frac 1k \, \min\{\ExpectedMisses,2\ExpectedMisses-2\mu\}$.
\end{lem}

We note that sometimes this bound is weaker than $z({|\Definite|+k} \mid \beta)\leq 1$. This happens for few enough $\beta$ that, from a computational vantage point, it is not worth checking for.

\begin{proof}
Trivially, $0 \leq z({|\Definite|+k} \mid \beta)$. We have
   \be \ExpectedMisses\ =\ \sum_{i=0}^\infty i\, \cdot \, z(|\Definite|+i) \geq k z(|\Definite|+k),\ee
whence $z(|\Definite|+k) \leq \ExpectedMisses/k$. But also,
   \begin{align}
      \ExpectedMisses  \ = \ \sum_{i=0}^\infty i\, \cdot \, z(|\Definite|+i) \ =\ z(|\Definite|+1) + \sum_{i=2}^\infty i\,\cdot\, z(|\Definite|+i) \ \geq\ 2\mu-\ExpectedMisses + k z(|\Definite|+k),
   \end{align}
and so $z(|\Definite|+k) \geq (2\ExpectedMisses-2\mu)/k$.
\end{proof}

%%%%%%%%%%%%%%%%%%%%%%%%%%%%%%%%%%%%%%%%%%%%%%%%%%%%%%%%%%%%%%%%%%%%%%%%%%%%%%%%%%%%%%%%%%%%%%%%%%%%%%%%%%%%%%%
\subsection{Making the computation feasible, reliable, and verifiable}

\newcommand{\Lower}{\ensuremath{\text{\sc Lower}}}
\newcommand{\Upper}{\ensuremath{\text{\sc Upper}}}
\newcommand{\Overhang}{\ensuremath{\text{\sc Overhang}}}

A massive computation has been performed, so some words are necessary as to how this is feasible. Set
   \begin{align}
     \Lower(k \mid \beta) &\ :=\
      \begin{cases}
         \max\{0,2^n-2^n\ExpectedMisses\}, & k = |\Definite| \\
         \max\{0,2\cdot2^n\mu-2^n\ExpectedMisses\}, & k = |\Definite|+1 \\
         0, & \text{otherwise}
      \end{cases} \nonumber\\
     \Upper(k \mid \beta) &\ :=\
      \begin{cases}
         2^n-2^n\mu, & k = |\Definite| \\
         \min\{2^n,2^n\ExpectedMisses\}, & k = |\Definite|+1 \\
         0, & \text{otherwise}
      \end{cases}\nonumber\\
     \Overhang(k | \beta) &\ :=\
      \begin{cases}
         \min\{2^n\ExpectedMisses,2\cdot 2^n\ExpectedMisses-2\cdot 2^n\mu\}, & k = |\Definite| \\
         0, & \text{otherwise}.
      \end{cases}
   \end{align}
The lemmas above imply that that the vector
   \be
      2^{2n-1} \langle z(0), z(1), \ldots, z(31)\rangle\ =\ \sum_{0\in \beta \subseteq[0,n)} 2^n \langle z(0\mid\beta),z(1\mid\beta),\ldots,z(31\mid\beta)\rangle
   \ee
is bounded below componentwise by
   \begin{equation}\label{equ: Lower}
      \sum_{0\in\beta\subseteq[0,n)} \langle \Lower(0\mid\beta),\Lower(1\mid\beta),\ldots,\Lower(31\mid\beta) \rangle
   \end{equation}
and is bounded above componentwise by
   \begin{multline*}
      \sum_{0\in\beta\subseteq[0,n)} \bigg(\langle \Upper(0\mid\beta),\Upper(1\mid\beta),\dots,\Upper(31\mid\beta)\rangle + \\
               \langle \Overhang(0\mid\beta),\Overhang(1\mid\beta),\ldots, \Overhang(31\mid\beta)\rangle \cdot M\bigg),
   \end{multline*}
where $M$ is the $32\times32$ matrix whose $(i,j)$\textsuperscript{th} entry (running the indices from 0 to 31) is $\frac{1}{j-i}$ if $j\geq i+2$, and is 0 otherwise.
This allows us to compute an upper bound on $z(0),\dots,z(31)$ from
   \begin{equation}\label{equ: Upper}
      \sum_{0\in\beta\subseteq[0,n)} \langle \Upper(0\mid\beta),\Upper(1\mid\beta),\dots,\Upper(31\mid\beta)\rangle
   \end{equation}
and
   \begin{equation}\label{equ: Overhang}
      \sum_{0\in\beta\subseteq[0,n)} \langle \Overhang(0\mid\beta),\Overhang(1\mid\beta),\dots,\Overhang(31\mid\beta)\rangle.
   \end{equation}
Observe that \Lower, \Upper{\rm\ and\ }Overhang\ are always integral, as $2^n\mu$ and $2^n\ExpectedMisses$ are both integers; this means that we can compute \eqref{equ: Lower}, \eqref{equ: Upper} and \eqref{equ: Overhang} using only integer arithmetic.

We need to compute $\beta+\beta$ and $\beta\cap[0,k]$ (for various $k$) for each $\beta$. This work can be tremendously reduced by using a Gray code. That is, the subsets of $[1,n)$ can be enumerated in such a way that each set differs from its predecessor in only one element (either put in or taken out). By storing the representation function for $\beta+\beta$ (that is, the number of times each sum can be written as a sum of two elements of $\beta$), we can simply update the necessary computations instead of re-computing.

Unfortunately, the size of the computation requires us to use $2n+1$-bit integers, and this is not a supported data type in most languages for $n\geq 32$. The options of using C with GMP, Mathematica, or some other route to arbitrary size integers is prohibited by the size of our computation and the modesty of our actual needs (we add, but never multiply, and know a priori the number of bits we will need). Therefore, we choose to represent our numbers as arrays of 64-bit integers in C++ (each element of the array represents a separate digit of the binary expansion of the number, but the digits aren't restricted to $\{0,1\}$). To further extend our reach, we ran the code on the parallel computing cluster at the High Performance Computing Cluster at the City University of New York. To facilitate parallelization, we break $\beta$ into $\beta_1=\beta\cap[0,n_1)$ and $\beta_2=\beta\cap[n_1,n)$. This makes the algorithm ``embarrassingly parallel'', and allows us to store intermediate calculations both to recover from any system or power failings, and to allow for spot checking of results.

To ensure correctness of the results, we have written the code in Mathematica using the simplest algorithms conceivable. Such code becomes intractably slow around $n\approx 25$, but this provides a sequence of values against which we can test our progressively more subtly written code, both in Mathematica and in C++. Our most sophisticated code is in C++.

Finally, we have the bounds on $\Prob{|\NN \setminus (D+D)|=k \mid D\cap[0,2^{10}) = \beta_1}$ for all $\beta_1$ in a publicly available file, together with our source code. We invite the reader to spot check our implementation.

%%%%%%%%%%%%%%%%%%%%%%%%%%%%%%%%%%%%%%%%%%%%%%%%%%%%%%%%%%%%%%%%%%%%%%%%%%%%%%%%%%%%%%%%%%%%%%%%%%%%%%%%%%%%%%%
\subsection{Obtaining \boldmath{$y(k)$}, \boldmath{$m(k)$}, and \boldmath{$w(k)$} from \boldmath{$z(k)$}}\label{subsec: connecting the distribs}

While it is clear that $z(k)$ is defined, that is, the event ``$|\NN\setminus(D+D)|=k$'' is measurable, it is less clear that $z(\infty)=0$. This, and that $y(\infty)=0$, follows from the Borel-Cantelli lemma and bounds such as~\eqref{eqn: iapp}. We can define $D$ (a uniformly chosen subset of $\NN$ containing 0) as $C-\min C$ (where $C$ is a uniformly chosen subset of $\NN$), and so
   \begin{align}
   y(k) &\ =\ \sum_{i=0}^\infty \Prob{\min C = i \text{ AND } |\NN\setminus(C+C)|=k} \nonumber\\
         &\ =\ \sum_{i=0}^\infty \Prob{\min C = i \text{ AND } |\NN\setminus((C-\min C)+(C-\min C))|=k-2i} \nonumber\\
         &\ =\ \sum_{i=0}^{\lfloor k/2 \rfloor} \Prob{\min C = i}\Prob{|\NN \setminus (D+D)| = k-2i} \ =\ \sum_{i=0}^{\lfloor k/2 \rfloor} \frac{1}{2^{i+1}} z(k-2i).
   \end{align}
To obtain the formulas
   \be
      m(k)\ =\ \sum_{i=0}^k y(i)y(k-i),\qquad \ \ w(k)\ =\ \sum_{i=0}^k z(i) z(k-i)
   \ee
we refer the reader to~\cite{In}. The gist of the argument is that
   \begin{align}
      m(k) &\ := \ \Prob{|[0,2n-2]\setminus(A+A)| = k}\nonumber\\
         &\ = \ \sum_{i=0}^k \Prob{|[0,n/2)\setminus(A+A)|=i\text{ AND } |(3n/2,2n-2] \setminus A+A| = k-i} \nonumber\\
         &\ \: \: \: \: \ \: \: \: \:  + \Prob{A+A \mbox{ misses an element in } [n/2, 3n/2]} \nonumber\\
          &\ = \ \sum_{i=0}^k \Prob{|[0,n/2)\setminus(A+A)|=i\text{ AND } |(3n/2,2n-2] \setminus A+A| = k-i}
          +O\left(\left(\frac{3}{4}\right)^{n/4} \right).
   \end{align}
Since $A+A\cap[0,n/2)$ is only affected by $A\cap[0,n/2)$ and $A+A\cap(3n/2,2n-2]$ is only affected by $A\cap(n/2,n)$, we can use independence to write
   \begin{equation}\label{ineq: appendix}
      m(k)\ = \ \sum_{i=0}^k \Prob{|[0,n/2)\setminus(A+A)|=i}\, \Prob{|(3n/2,2n-2] \setminus A+A| = k-i} + O\left(\left(\frac{3}{4}\right)^{n/4} \right).
   \end{equation}
so that
   \be
      m(k)\ \sim \ \sum_{i=0}^k \Prob{|[0,n/2)\setminus(A+A)|=i}\, \Prob{|(3n/2,2n-2] \setminus A+A| = k-i}.
   \ee
As $n\to\infty$, the set $[0,n/2)\setminus (A+A)$ looks more and more like $\NN \setminus (C+C)$, so that
   \be
      \Prob{|[0,n/2)\setminus(A+A)|=i} \to y(i),
   \ee
and similarly (after replacing $A$ with $n-1-A$) for $ \Prob{|(3n/2,2n-2] \setminus A+A| = k-i}$. The argument for $w(k)$ is identical, but with ``$D$'' in place ``$C$''.

Let $Z_1,Z_2$ be independent random variables with the same distribution as $M_{\NN \mid \{0\}}$, and set $W:=Z_1+Z_2$. Then $\Prob{W=k}=\sum_{i=0}^k z(i)z(k-i)= w(k)$, whence $\sum_{i=0}^\infty w(i) =1 $, and similarly $\sum_{i=0}^\infty m(i) = 1$.

Since $y(k)$ is a linear combination of $z(0),\dots,z(k)$ with {\em positive} coefficients, the lower bounds on $z(0)$, $\dots$, $z(k)$ immediately give a lower bound on $y(k)$, and likewise upper bounds on $z(0)$, $\dots$, $z(k)$ yield an upper bound on $y(k)$. The situation is the same between $y$ and $m$ and between $z$ and $w$, even though the combination is not linear!

To experimentally estimate $z(k)$, we hypothesized that $\Prob{\NN \setminus(D+D) \not\subseteq [0,256)}$ is sufficiently small as to be ignored. Then, using Mathematica 8, we generated $2^{28}$ pseudorandom subsets $E$ of $[0,256)$, forced each to contain 0, and then computed $k:=|[0,256) \setminus (E+E)|$ and kept a running tally of the number of times each value of $k$ arose. This estimates (with an enormous sample size) \be \Prob{|\NN \setminus (D+D)| = k \mid \NN \setminus(D+D)\ \subseteq\ [0,256)}\approx z(k).\ee
The estimates $\widehat{z(k)}$, along with conservative 99.9\% confidence intervals, are given in Table~\ref{table for z} and shown in Figure~\ref{figure:z_k Experimental graph}. The implied bounds on $w$, $m$, and $y$ are given in Tables~\ref{table for y},~\ref{table for m}, and~\ref{table for w} respectively, and shown in Figure~\ref{figure:z_k Experimental graph}.

%%%%%%%%%%%%%%%%%%%%%%%%%%%%%%%%%%%%%%%%%%%%%%%%%%%%%%%%%%%%%%%%%%%%%%%%%%%%%%%%%%%%%%%%%%%%%%%%%%%%%%%%%%%%%%%

%%%%%%%%%%%%%%%%%%%%%%%%%%%%%%%%%%%%%%%%%%%%%%%%%%%%%%%%%%%%%%%%%%%%%%%%%%%%%%%%%%%%%%%%%%%%%%
%%%%%%%%%%%%%%%%%%%%%%%%%%%%%%%%%%%%%%%%%%%%%%%%%%%%%%%%%%%%%%%%%%%%%%%%%%%%%%%%%%%%%%%%%%%%%%
\section{Conjectures and Future Research}\label{sec: future}

We end with some conjectures that are supported by numerical data. Our main conjecture remains Conjecture \ref{conj: main}, which says that the distribution of missing sums is approximately exponential.
One possible method of studying this distribution is finding where the first present sum in $A+A$ occurs, given that $A+A$ has $k$ missing elements.
Recall that the lower bound in \S\ref{sec: bounds} was proven by constructing $A$ such that $M_{[0, n-1]}(A) = k$ by letting the first $k/2$ elements of $A$ be missing.
In this case, the index of the first present sum in $A+A$ occurs at index $k$.
But from numerical data, the index of the first present element will not be $k$ for typical $A+A$ that is missing $k$ elements.
This also suggests that this trivial construction does not account for the real `random' way of constructing $A$ such that $A+A$ is missing $k$ elements, which is consistent with the fact that the conjectured decay constant for the distribution is approximately $0.78$ but the lower bound gives only the decay constant approximately $0.70$. Even though the index of the first present element is not $k$, from numerical data, the index seems to be linear in $k$.

To be precise, let $$X_n(A)\  =\ \max \{m: {\rm if}\ \ell < m \ {\rm then}\ \ell \not\in A+A \}$$ be the index of the first present sum of $A+A$. Then we have the following conjecture.

\begin{conj}\label{conj1}
For large $k$,
\begin{equation}
\lim_{n\rightarrow \infty} \E( X_n(A) \mid M_{[0, n-1]}(A) =k )
\end{equation}
is asymptotically linear in $k$.
\end{conj}

%Here we have a limit in $n$ largely to emphasize that the terms do not depend on $n$ but only on $k$;
%however, it may be possible to prove that for any fixed $k$, this limit exists. One possible restatement of the conjecture is the following:
%there exist $a,b$ such that
%\begin{equation}
% \E( X_n(A) \mid M_{[0, n-1]}(A) = k) \sim ax+b,
%\end{equation}
%where $n > C k$ for some constant $C$.

Similarly, we can investigate how far we must move to the right of zero and to the left of the maximum possible sum, $2n-2$, so that there are no missing sums of $A+A$ in this interval. Given $A \in [0, n-1]$ missing exactly $k$ sums, as $n\to\infty$ each of the $k$ missing elements of $A+A$ are either near 0 or near $2n-2$. Thus all of the action is happening near the two fringes, and we want to understand what is happening there. This suggests studying $$\max \left\{Y_n(A) - W_n(A): [W_n(A), Y_n(A)] \subset A+A\right\}.$$

%consider the first index by which $A+A$ has all of its $k$ missing elements; that is,
%$Y_n(A) = \min\{m:  [0, m-1] \cap (A+A) \mbox{ has $k$ missing sums}\}$.
%If $Y_n(A) = k$, then $A+A$ is missing exactly the first $k$ elements, which is the same situation as in the trivial lower bound; however, for %typical $A$ this does not seem to be the case, and we conjecture

\begin{conj} With $W_n(A)$ as above 
\begin{equation}
\lim_{n\rightarrow \infty} \E( W_n(A) \mid M_{[0, n-1]}(A) =k )
\end{equation}
is asymptotically linear in $k$. 
\end{conj}

Note a similar conjecture should hold for $2n-2- Y_n(A)$.

%As for Conjecture \ref{conj1}, this conjecture has the following restatement:
%there exist $a,b$ such that
%\begin{equation}
% \E( Y_n(A) \mid M_{[0, n-1]}(A) = k) \sim ax+b,
%\end{equation}
%where $n > C k$ for some constant $C$.

Another direction is to improve the exponential bounds for  $\pp(M_{[0, n-1]}(A)= k)$. One approach to do this is to find upper bounds on probabilities like $\pp(a_1, \dots, a_m\not \in A+A)$ for arbitrary $a_1,a_2,\dots, a_m$ around $k$.

Recall that in \S\ref{sec: bounds} we first used
$\pp(i \not \in A+A)$ to get an upper bound for
$\pp(M_{[0, n-1]}(A)= k)$ of $\Theta \left( (3/4)^{k/2} \right)$ and
then used $\pp(i,j \not \in A+A)$ to get a bound of $\Theta \left( (\phi/2)^{k} \right)$, an improvement.
Knowing $\pp(a_1, \dots, a_m\not \in A+A)$ would result in similar improvement. Using the current approach, this would require studying the number of vertex covers for graphs that have vertices with degree $m$ instead of $2$.

\begin{comment}
Note that we already have an upper bound $\pp(k+1,\dots, k+m\not \in A+A)$ from \S\ref{sec: consecutive} since this is the case of consecutive missing sums.
This probability tends towards $(1/2)^{(k+m)/2}$ (in the sense of Theorem \ref{thm: consecutive}), which matches the trivial lower bound.
If  $\pp(k+a_1, k+a_2,\dots, k+a_m\not \in A+A) \le \pp(k+1,\dots, k+m \not \in A+A)$ for any $a_1, \dots, a_m$, then the probability for missing $k$ elements would be $(1/2)^{(k+m)/2}$. Since this does not seem to be the case from numerical data, the probability of missing $m$ elements around $k$ is asymptotically different from the probability of missing $m$ consecutive elements around $k$; that is $\lambda_{a_1, \dots, a_m}$ from Theorem \ref{thm: configuration}
is not close to $(1/2)^{1/2}$.
\end{comment}

Finally, note that it is possible to use the graph-theoretic approach to study higher moments of $M_{[0, n-1]}$. Recall that the variance was calculated by finding explicit formulas for $\pp(i \mbox{ and } j \not \in A+A)$. Similarly, the $m$\textsuperscript{th} moment can be found by finding explicit formulas for
$\pp(a_1, \dots, a_m \not \in A+A)$ for arbitrary $a_1, \dots, a_m$, which requires finding the number of vertex covers in certain graphs that have vertices with degree $m$.
Note that we again need to study $\pp(a_1, \dots, a_m \not \in A+A)$, as we do when we try to improve the bounds for $\pp(M_{[0, n-1]}(A) = k)$;
however now we need an exact formula for $\pp(a_1, \dots, a_m \not \in A+A)$, whereas before we just needed an upper bound.

%%%%%%%%%%%%%%%%%%%%%%%%%%%%%%%%%%%%%%%%%%%%%%%%%%%%%%%%%%%%%%%%%%%%%%%%%%%%%%%%%%%%%%%%%%%%%%

%%%%%%%%%%%%%%%%%%%%%%%%%%%%%%%%%%%%%%%%%%%%%%%%%%%%%%%%%%%%%%%%%%%%%%%%%%%%%%%%%%%%%%%%%%%%%%

\newpage

\appendix

\section{Data tables for distributions}\label{Appendix}

\begin{figure}[h]
   \[
      \begin{array}{c|ccccc}
        & \text{rigorous} & \text{lower} &  & \text{upper} & \text{rigorous} \\
      k & \text{lower} & \text{CI}  & 10^5\widehat{z(k)} & \text{CI} & \text{upper} \\ \hline
 0 & 23532 & 23543 & 23554 & 23566 & 23535 \\
 1 & 17651 & 17634 & 17644 & 17655 & 17662 \\
 2 & 13955 & 13941 & 13950 & 13960 & 13975 \\
 3 & 11074 & 11065 & 11073 & 11082 & 11101 \\
 4 & 9233 & 9225 & 9233 & 9241 & 9266 \\
 5 & 6502 & 6502 & 6509 & 6516 & 6540 \\
 6 & 5049 & 5055 & 5061 & 5067 & 5090 \\
 7 & 3700 & 3710 & 3716 & 3721 & 3745 \\
 8 & 2687 & 2698 & 2703 & 2708 & 2733 \\
 9 & 1898 & 1910 & 1914 & 1918 & 1945 \\
 10 & 1384 & 1400 & 1404 & 1407 & 1433 \\
 11 & 958 & 973 & 976 & 979 & 1006 \\
 12 & 677 & 691 & 694 & 697 & 725 \\
 13 & 467 & 480 & 483 & 485 & 515 \\
 14 & 323 & 337 & 339 & 341 & 370 \\
 15 & 219 & 231 & 233 & 235 & 266 \\
 16 & 149 & 161 & 162 & 164 & 195 \\
 17 & 100 & 110 & 111 & 112 & 145 \\
 18 & 66 & 75 & 76 & 77 & 110 \\
 19 & 43 & 51 & 52 & 53 & 86 \\
 20 & 28 & 35 & 36 & 37 & 70 \\
 21 & 18 & 23 & 24 & 25 & 58 \\
 22 & 11 & 16 & 16 & 17 & 51 \\
 23 & 7 & 11 & 11 & 12 & 45 \\
 24 & 4 & 7 & 8 & 8 & 42 \\
 25 & 2 & 4 & 5 & 6 & 39 \\
 26 & 1 & 3 & 4 & 4 & 37 \\
 27 & 0 & 2 & 2 & 3 & 36 \\
 28 & 0 & 1 & 2 & 2 & 35 \\
 29 & 0 & 1 & 1 & 2 & 35 \\
 30 & 0 & 0 & 1 & 1 & 34 \\
 31 & 0 & 0 & 1 & 1 & 34 \\
      \end{array}
   \]
\caption{The first and last columns give our rigorous lower and upper bounds on $10^5 z(k)$. The second and fourth columns give the bounds of a conservative 99.9\% confidence interval for $10^5\widehat{z(k)}$. The middle column gives our best guess for the integer closest to $10^5z(k)$, which we denote $10^5 \widehat{z(k)}$.
\label{table for z}}
\end{figure}

\begin{figure}[h]
   \[
      \begin{array}{c|ccccc}
            & \text{rigorous} & \text{lower} &  & \text{upper} & \text{rigorous} \\
         k & \text{lower} & \text{CI}  & 10^5\widehat{y(k)} & \text{CI} & \text{upper} \\ \hline
 0 & 11766 & 11771 & 11777 & 11783 & 11768 \\
 1 & 8825 & 8817 & 8822 & 8828 & 8831 \\
 2 & 12860 & 12856 & 12864 & 12871 & 12872 \\
 3 & 9950 & 9941 & 9948 & 9955 & 9966 \\
 4 & 11047 & 11041 & 11048 & 11056 & 11069 \\
 5 & 8226 & 8221 & 8228 & 8235 & 8253 \\
 6 & 8048 & 8048 & 8055 & 8062 & 8079 \\
 7 & 5963 & 5966 & 5972 & 5978 & 5999 \\
 8 & 5367 & 5373 & 5379 & 5385 & 5406 \\
 9 & 3931 & 3938 & 3943 & 3948 & 3972 \\
 10 & 3376 & 3387 & 3391 & 3396 & 3419 \\
 11 & 2444 & 2455 & 2459 & 2463 & 2489 \\
 12 & 2026 & 2039 & 2043 & 2046 & 2072 \\
 13 & 1456 & 1468 & 1471 & 1474 & 1502 \\
 14 & 1174 & 1188 & 1191 & 1193 & 1221 \\
 15 & 837 & 850 & 852 & 855 & 884 \\
 16 & 662 & 674 & 676 & 679 & 708 \\
 17 & 468 & 480 & 482 & 483 & 514 \\
 18 & 364 & 375 & 376 & 378 & 409 \\
 19 & 256 & 265 & 267 & 268 & 300 \\
 20 & 196 & 205 & 206 & 207 & 240 \\
 21 & 137 & 144 & 146 & 147 & 179 \\
 22 & 103 & 110 & 111 & 112 & 145 \\
 23 & 72 & 77 & 78 & 79 & 112 \\
 24 & 54 & 59 & 59 & 60 & 93 \\
 25 & 37 & 41 & 42 & 43 & 76 \\
 26 & 27 & 31 & 32 & 32 & 65 \\
 27 & 19 & 21 & 22 & 23 & 56 \\
 28 & 14 & 16 & 17 & 17 & 50 \\
 29 & 9 & 11 & 12 & 12 & 45 \\
 30 & 7 & 8 & 9 & 9 & 42 \\
 31 & 4 & 5 & 6 & 7 & 40 \\
     \end{array}
   \]
\caption{The first and last columns give our rigorous lower and upper bounds on $10^5 y(k)$. The second and fourth columns give the bounds of a conservative 99.9\% confidence interval for $10^5\widehat{y(k)}$. The middle column gives our best guess for the integer closest to $10^5y(k)$, which we denote $10^5 \widehat{y(k)}$.
\label{table for y}}
\end{figure}

\begin{figure}[h]
   \[
      \begin{array}{c|ccccc}
        & \text{rigorous} & \text{lower} &  & \text{upper} & \text{rigorous} \\
      k & \text{lower} & \text{CI}  & 10^5\widehat{m(k)} & \text{CI} & \text{upper} \\ \hline
 0 & 1384 & 1385 & 1387 & 1389 & 1385 \\
 1 & 2076 & 2075 & 2078 & 2081 & 2079 \\
 2 & 3805 & 3804 & 3808 & 3813 & 3810 \\
 3 & 4611 & 4607 & 4613 & 4618 & 4619 \\
 4 & 6010 & 6005 & 6012 & 6020 & 6022 \\
 5 & 6445 & 6439 & 6447 & 6455 & 6463 \\
 6 & 7177 & 7172 & 7181 & 7191 & 7202 \\
 7 & 7138 & 7133 & 7143 & 7153 & 7170 \\
 8 & 7243 & 7240 & 7251 & 7261 & 7282 \\
 9 & 6825 & 6824 & 6835 & 6846 & 6871 \\
 10 & 6510 & 6513 & 6523 & 6534 & 6563 \\
 11 & 5892 & 5897 & 5907 & 5918 & 5951 \\
 12 & 5374 & 5382 & 5392 & 5402 & 5439 \\
 13 & 4712 & 4724 & 4733 & 4742 & 4783 \\
 14 & 4153 & 4168 & 4176 & 4185 & 4228 \\
 15 & 3551 & 3567 & 3575 & 3583 & 3629 \\
 16 & 3046 & 3064 & 3071 & 3079 & 3127 \\
 17 & 2550 & 2569 & 2576 & 2582 & 2633 \\
 18 & 2139 & 2159 & 2165 & 2172 & 2225 \\
 19 & 1759 & 1780 & 1785 & 1790 & 1846 \\
 20 & 1449 & 1469 & 1474 & 1479 & 1536 \\
 21 & 1173 & 1193 & 1198 & 1202 & 1260 \\
 22 & 951 & 970 & 974 & 978 & 1038 \\
 23 & 760 & 778 & 782 & 785 & 846 \\
 24 & 608 & 625 & 628 & 631 & 693 \\
 25 & 480 & 496 & 498 & 501 & 564 \\
 26 & 379 & 394 & 396 & 398 & 462 \\
 27 & 296 & 309 & 311 & 313 & 378 \\
 28 & 232 & 243 & 245 & 247 & 311 \\
 29 & 179 & 189 & 191 & 193 & 258 \\
 30 & 139 & 148 & 149 & 150 & 216 \\
 31 & 106 & 114 & 115 & 117 & 182 \\
      \end{array}
   \]
\caption{The first and last columns give our rigorous lower and upper bounds on $10^5 m(k)$. The second and fourth columns give the bounds of a conservative 99.9\% confidence interval for $10^5\widehat{m(k)}$. The middle column gives our best guess for the integer closest to $10^5m(k)$, which we denote $10^5 \widehat{m(k)}$.
\label{table for m}}
\end{figure}

\begin{figure}[h]
   \[
      \begin{array}{c|ccccc}
        & \text{rigorous} & \text{lower} &  & \text{upper} & \text{rigorous} \\
      k & \text{lower} & \text{CI}  & 10^5\widehat{w(k)} & \text{CI} & \text{upper} \\ \hline
 0 & 5537 & 5543 & 5548 & 5554 & 5539 \\
 1 & 8307 & 8303 & 8312 & 8321 & 8314 \\
 2 & 9684 & 9674 & 9685 & 9696 & 9698 \\
 3 & 10138 & 10127 & 10139 & 10152 & 10162 \\
 4 & 10202 & 10190 & 10203 & 10217 & 10236 \\
 5 & 9411 & 9401 & 9414 & 9427 & 9454 \\
 6 & 8475 & 8470 & 8483 & 8497 & 8528 \\
 7 & 7384 & 7385 & 7397 & 7410 & 7445 \\
 8 & 6273 & 6279 & 6291 & 6302 & 6342 \\
 9 & 5194 & 5204 & 5215 & 5226 & 5269 \\
 10 & 4247 & 4262 & 4272 & 4282 & 4327 \\
 11 & 3405 & 3424 & 3433 & 3441 & 3490 \\
 12 & 2696 & 2718 & 2726 & 2733 & 2784 \\
 13 & 2107 & 2130 & 2137 & 2144 & 2197 \\
 14 & 1629 & 1654 & 1660 & 1666 & 1720 \\
 15 & 1245 & 1270 & 1275 & 1281 & 1337 \\
 16 & 943 & 968 & 973 & 977 & 1035 \\
 17 & 708 & 732 & 736 & 740 & 800 \\
 18 & 527 & 549 & 553 & 556 & 617 \\
 19 & 389 & 410 & 412 & 415 & 478 \\
 20 & 285 & 304 & 306 & 309 & 372 \\
 21 & 207 & 224 & 226 & 228 & 293 \\
 22 & 149 & 164 & 166 & 168 & 233 \\
 23 & 106 & 120 & 121 & 123 & 189 \\
 24 & 75 & 87 & 88 & 90 & 156 \\
 25 & 53 & 63 & 64 & 65 & 132 \\
 26 & 37 & 45 & 46 & 48 & 114 \\
 27 & 25 & 32 & 33 & 34 & 101 \\
 28 & 17 & 23 & 24 & 25 & 91 \\
 29 & 12 & 16 & 17 & 18 & 84 \\
 30 & 8 & 11 & 12 & 13 & 79 \\
 31 & 5 & 8 & 9 & 10 & 76 \\
      \end{array}
   \]
\caption{The first and last columns give our rigorous lower and upper bounds on $10^5 w(k)$. The second and fourth columns give the bounds of a conservative 99.9\% confidence interval for $10^5\widehat{w(k)}$. The middle column gives our best guess for the integer closest to $10^5w(k)$, which we denote $10^5 \widehat{w(k)}$.
\label{table for w}}
\end{figure}

%%%%%%%%%%%%%%%%%%%%%%%%%%%%%%%%%%%%%%%%%%%%%%%%%%%%%%%%%%%%%%%%%%%%%%%%%%%%%%%%%%%%%%%%%%%%%%
%%%%%%%%%%%%%%%%%%%%%%%%%%%%%%%%%%%%%%%%%%%%%%%%%%%%%%%%%%%%%%%%%%%%%%%%%%%%%%%%%%%%%%%%%%%%%%

\ \\

\end{document}